\documentclass[12pt,a4paper]{article}
\usepackage[cp1251]{inputenc}
\usepackage[T2A]{fontenc}
\usepackage[english]{babel}

\usepackage{amsmath,amsfonts,amssymb,mathrsfs,amsopn,indentfirst,amscd}

\usepackage{graphicx}
\usepackage{amsthm}
\usepackage{hyperref}

\usepackage{amssymb}
\usepackage{amsmath}

\newcommand{\ver}{{\rm ver}}
\newcommand{\vol}{{\rm vol}}
\newcommand{\vo}{{\rm vol}}

\newcommand{\mes}{{\rm mes_n}}

\newtheorem*{corollary*}{Corollary}

\begin{document}

\title{Estimates for Interpolation Projectors \\
and Related Problems \\  in Computational Geometry\\}
\author{Mikhail Nevskii\footnote{ Department of Mathematics,  P.\,G.~Demidov Yaroslavl State University, Sovetskaya str., 14, Yaroslavl, 150003, Russia, 
               mnevsk55@yandex.ru,  orcid.org/0000-0002-6392-7618 } 
               \ and
               Alexey Ukhalov\footnote{ Department of Mathematics,  P.\,G.~Demidov Yaroslavl State University, Sovetskaya str., 14, Yaroslavl, 150003, Russia, 
               alex-uhalov@yandex.ru,  orcid.org/0000-0001-6551-5118 }  
               }
      
\date{August 2, 2021}
\maketitle


\begin{quotation}

This paper contains a survey of  results obtained by the authors mostly during the past few years and published by 2021.
In particular, we present   the best of  known estimates of numerical 
characteristics related  to the research theme.


\smallskip
{\bf Sections:} 1. Introduction. 2. The case when $n+1$ is an Hadamard number. 3. Estimates for the minimal absorption index of a cube by a~simplex. 4.~Estimates for the minimal norm of a projector
 in linear\-  interpolation on a~cube in ${\mathbb R}^n$. 5. Estimates of numbers $\xi_n^\prime$ and $\theta_n^\prime$.
6.~Simplices satisfying the~inclusions $S\subset Q_n\subset nS$. 7. Perfect simplices. 8.~Equisecting 
 simplices.
9. Properties of $(0,1)$-matrices of order $n$ having maximal determinant. 10.~Problems for a simplex\linebreak  and 
a Euclidean ball. 11. Linear interpolation on a Euclidean ball.

\smallskip
{\bf Bibliography:} 56 titles.



\end{quotation}

\section{Introduction}\label{nev_ukh_sec_1}

This paper provides a survey  of some results obtained by the authors mainly in the last seven years.
The subject of our consideration is
polynomial interpolation in several variables and also application in this area of geometric constructions and methods. As analysis in general,
approximation theory is closely related to geometry. It is well known that many fundamental
 results in approximation theory 
essentially have geometric character. In detail,  this connection   is discovered  in  literature.

In polynomial interpolation of functions in several variables, the important geometric questions arise in connection with the structure of the basic domain and the choice of the nodes set.
In particular, when interpolating
functions of $n$ variables using the space of
polynomials of degree $ \leq 1 $, 
the interpolation nodes coincide with
vertices of an $n$-dimensional nondegenerate simplex. By this, it occurs to be possible to provide estimates of the interpolation projector's norm through some geometric or numerical characteristics of the corresponding simplex.
This is where we can use the properties of some numeric characteristics of a simplex and another convex body
(volume, axial diameters, minimal positive homothety coefficients related to absorption  a convex body  by a homothet of a simplex,  with or without translation).
This approach can be carried over also to interpolation 
with  wider spaces of~polynomials.

The basic theoretical apparatus in this area was developed by M.\,Nevsky approximately
in 2000--2014. Here we  note  the monograph
\cite{nevskii_monograph} containing many results and links, the  papers
\cite{nevskii_mais_2003_10_1}--\cite{nevskii_mais_2006_13_2},
\cite{nevskii_mais_2007_14_1}--\cite{nevskii_dcg_2011},
\cite{nevskii_matzam_2011},
\cite{nevskii_matzam_2013}--\cite{nevskii_fpm_2013},
and also  the dissertation~\cite{nevskii_disser}. 
A.\,Ukhalov  joined the research
 starting around 2015. In particular, he very successfully applied computer methods in this area.
 The is written mainly on the basis of papers published by the authors or their students in 2016--2020. 
In order to make the text  more understandable, 
we also cite the works written in previous years.

This paper continues the series of review articles in the anniversary collections issued 
for the 25th, 30th,  35th, and 40th anniversaries of Department of Mathematics of  P.\,G.\,Demidov Yaroslavl State University (see  
 \cite{irodova_nevskii_sb_2001},  \cite{nevskii_sb_2006},   \cite{nevskii_sb_2011},  \cite{nevskii_sb_2016}).
 The   original Russian text was prepared for a collection dedicated to 45th anniversary of the  faculty which  took place in 2021.

\bigskip
We start with main definitions and denotations.  Everywhere further $n\in{\mathbb N}.$ An element
$x\in{\mathbb R}^n$ is written in the form
$x=(x_1,\ldots,x_n).$ 
By $e_1$, $\ldots$, $e_n$ we denote the standard basis
in ${\mathbb R}^n$.
By definition,  $$\|x\|:=\sqrt{(x,x)}=\left(\sum\limits_{i=1}^n x_i^2\right)^{1/2},$$ 
$$B\left(x^{(0)};R\right):=\{x\in{\mathbb R}^n: \|x-x^{(0)}\|\leq R \} 
\quad \left(x^{(0)}\in {\mathbb R}^n,
R>0\right),$$ \ 
$$B_n:=B(0;1), \quad 
Q_n:=[0,1]^n, \quad
Q_n^\prime:=[-1,1]^n.$$
The notation 
$L(n)\asymp M(n)$ means that there exist  $c_1,c_2>0$
not depending on $n$ such that
$c_1M(n)\leq L(n)\leq c_2 M(n)$.

Let $\Omega$ be   {\it a convex body in ${\mathbb R}^n,$ } i.\,e., a compact convex subset of ${\mathbb R}^n$ 
with nonempty interior. 
By $\sigma \Omega$ we denote a homothetic copy of $C$ with the center of homothety in the center of
gravity of $C$ and the ratio of homothety  $\sigma.$ 
The symbol $\vol(\Omega)$ denotes the volume of
$\Omega.$
If $\Omega$ is a convex polytope, then $\ver(\Omega)$ is a set
of vertices of  $\Omega$. By {\it a translate} we mean the result of a parallel shift.

An $n$-dimensional simplex $S$ is circumscribed around a convex body
$\Omega$ if $\Omega\subset S$ and each
$(n-1)$-dimensional face of  $S$ contains a point of $\Omega.$ 
A convex polytope is inscribed into $\Omega$ if every vertex of this polytope belongs to the boundary of
$\Omega.$

 Define 
$d_i(\Omega)$ as the minimal length of a segment contained in $\Omega$ and parallel to the 
$x_i$-axis. 
We call 
$d_i(\Omega)$  {\it the $i$th
axial diameter of $\Omega$}. 
The notion {\it axial diameter of a convex body} was introduced by P.\,Scott \cite{scott_1985},\cite{scott_1989}. 

For convex bodies
$\Omega_1$,
$\Omega_2$, let us denote by  
$\xi(\Omega_1;\Omega_2)$ 
 the minimal  $\sigma\geq 1$ having the property
$\Omega_1\subset \sigma \Omega_2$. An equality 
$\xi(\Omega_1;\Omega_2)$ $=1$ is  equivalent to the inclusion  \linebreak$\Omega_1\subset \Omega_2$.
We call
$\xi(\Omega_1,\Omega_2)$ {\it the absorption index of $\Omega_1$
by  $\Omega_2$}. 
Define 
$\alpha(\Omega_1,\Omega_2)$ as minimal $\sigma>0$ such that
$\Omega_1$ is a subset of a translate of
 $\sigma \Omega_2$. 
Clearly, $\alpha(\Omega_1,\Omega_2)$ $\leq $ $\xi(\Omega_1,\Omega_2)$.We put
$\xi(\Omega):=\xi(Q_n;\Omega),$
$\alpha(\Omega):=\alpha(Q_n;\Omega)$. 

The symbol
$C(\Omega)$ denotes a space of continuous functions
$f:\Omega\to{\mathbb R}$ with the uniform norm
$$\|f\|_{C(\Omega)}:=\max\limits_{x\in \Omega}|f(x)|.$$
By $\Pi_1\left({\mathbb R}^n\right)$ we mean a set of polynomials in 
$n$ variables of degree
$\leq1$, i.\,e., a set of linear functions upon  ${\mathbb R}^n$.

Let $S$ be a nondegenerate simplex in ${\mathbb R}^n$ with vertices             
$x^{(j)}=\left(x_1^{(j)},\ldots,x_n^{(j)}\right),$ 
$1\leq j\leq n+1.$  
{\it The vertices matrix}  of this simplex
$${\bf A} :=
\left( \begin{array}{cccc}
x_1^{(1)}&\ldots&x_n^{(1)}&1\\
x_1^{(2)}&\ldots&x_n^{(2)}&1\\
\vdots&\vdots&\vdots&\vdots\\
x_1^{(n+1)}&\ldots&x_n^{(n+1)}&1\\
\end{array}
\right)$$
is nondegenerate and $\vol(S)=\frac{1}{n!}|\det({\bf A})|.$  Assume 
${\bf A}^{-1}$ $=(l_{ij})$. Define  $\lambda_j\in  \Pi_1({\mathbb R}^n)$ as polynomials with
coefficients forming the columns of
${\bf A}^{-1}$:
$$\lambda_j(x)=
l_{1j}x_1+\ldots+
l_{nj}x_n+l_{n+1,j}.$$
These polynomials have the property  
$\lambda_j\left(x^{(k)}\right)$ $=$ 
$\delta_j^k$,
where $\delta_j^k$ is the Kronecker delta. 
We call $\lambda_j$ 
{\it the basic Lagrange polynomials corresponding to  $S$.}
For $x\in {\mathbb R}^n$, the  numbers $\lambda_1(x),$ $\ldots,$ $\lambda_{n+1}(x)$  
are 
{\it the barycentric coordinates of $x$} with respect to  
$S$.  
Equations $\lambda_j(x)=0$ give the  $(n-1)$-dimensional hyperplanes  containing the faces of
 $S.$  We have $$
S=\left \{ x\in {\mathbb R}^n: \, \lambda_j(x) \geq 0, \, j=1,\ldots,n+1
\right\}. 
$$
For more information on $\lambda_j$, see ~\cite{nevskii_monograph}, \cite{nev_ukh_posobie_2020}.

The first author proved that
\begin{equation}\label{xi_C_S_formula}
\xi(\Omega;S)=(n+1)\max_{1\leq k\leq n+1}
\max_{x\in \Omega}(-\lambda_k(x))+1 \quad (\Omega\not\subset S),   
\end{equation}
\begin{equation}\label{alpha_C_S_general_formula}
\alpha(\Omega;S)=\sum_{j=1}^{n+1} \max_{x\in \Omega} (-\lambda_j(x))+1.
\end{equation}
(see \cite{nevskii_dcg_2011}; the proofs are given also in \cite{nevskii_monograph}).
In the present text, we consider  mainly the cases $\Omega=Q_n$  and $\Omega=B_n$.

If  $\Omega=Q_n$, then \eqref{xi_C_S_formula}  and  \eqref{alpha_C_S_general_formula} are equivalent
correspondingly to the equalities
\begin{equation}\label{xi_S_formula}
\xi(S)=(n+1)\max_{1\leq k\leq n+1}
\max_{x\in \ver(Q_n)}(-\lambda_k(x))+1,
\end{equation}
\begin{equation}\label{alpha_S_formula}
\alpha(S)=
\sum_{j=1}^{n+1} \max_{x\in \ver(Q_n)} (-\lambda_j(x))+1.
\end{equation}
The equality $\xi(S)=\alpha(S)$  holds true if and only if
the simplex $\xi(S)S$ is circumscribed around $Q_n$. This  is also equivalent to the relation 
\begin{equation}\label{xi_S_S_circ_condition_around_Q_n}
 \max_{x\in \ver(Q_n)} \left(-\lambda_1(x)\right)=
\ldots=
\max_{x\in \ver(Q_n)} \left(-\lambda_{n+1}(x)\right).
\end{equation}

Very important result is the following formula for  $\alpha(S)$ obtained in \cite{nevskii_dcg_2011}:
\begin{equation}\label{alpha_d_i_formula}
\alpha(S)=\sum_{i=1}^n\frac{1}{d_i(S)}. 
\end{equation}
It was proved in  \cite{nevskii_matzam_2010} that
\begin{equation}\label{d_i_l_ij_formula}
\frac{1}{d_i(S)}=\frac{1}{2}\sum_{j=1}^{n+1} |l_{ij}|.
\end{equation}
From (\ref{alpha_d_i_formula}) and (\ref{d_i_l_ij_formula}), we have
\begin{equation}\label{alpha_qs_formula}
\alpha(S)=\frac{1}{2}\sum_{i=1}^n\sum_{j=1}^{n+1} |l_{ij}|. 
\end{equation}
Obviously, $\alpha(\Omega_1;\Omega_2)$ is invariant under the set translations and 
$\alpha(\tau \Omega_1;\Omega_2)$ $=$ $\tau\alpha(\Omega_1;\Omega_2)$, $\tau>0$.  Since 
$Q_n^\prime=[-1,1]^n$  
is a translate of $2Q_n$, equality
 (\ref{alpha_qs_formula}) 
gives
\begin{equation}\label{alpha_q_prime_s_formula}
\alpha(Q_n^\prime;S)=\sum_{i=1}^n\sum_{j=1}^{n+1} |l_{ij}|.
\end{equation}

Relation  \eqref{alpha_d_i_formula} has a lot of  useful corollaries given in
\cite{nevskii_dcg_2011} and \cite{nevskii_monograph}.
For instance, from  \eqref{alpha_d_i_formula}  it follows 
that for a simplex $S\subset Q_n$
\begin{equation}\label{xi_S_alpha_S_geq_n}
\xi(S)\geq \alpha(S)=\sum_{i=1}^n\frac{1}{d_i(S)}\geq n,
\end{equation}
since $d_i(S)\leq 1$. If
$\xi(S)=n$, then $\alpha(S)=n$ и $d_i(S)=1.$

Let us note also a result obtained on another way by M.\,Lassak \cite{lassak_dsg_1999}.
 {\it If a simplex $S\subset Q_n$ has the maximum volume, then $d_i(S)=1$.}
Indeed, for such an $S$ 
holds  $S\subset -nS$. If this was not so, then  some vertex of $S$ may be moved in $Q_n$ in such a way that
its distance from the opposite face of the simplex increases.
In this case the volume of $S$ was not maximum. The inclusion 
$Q_n\subset -nS$ means that  $-Q_n$ is a subset of a translate of
$nS$. Since  $-Q_n=Q_n$, we have $\alpha(S)\leq n$. As it was marked, if  $S\subset Q_n$, then $\alpha(S)\geq n$.
Therefore, $\alpha(S)=n$ and   
\eqref{alpha_d_i_formula} gives $d_i(S)=1.$ This property of a maximum volume simplex essentially was used in
\cite{nev_ukh_matzam SVFU_2019}.

Computational formullae for $\xi(B_n;S)$ and $\alpha(B_n;S)$ were obtained by M.\,Nevskii in  \cite{nevskii_mais_2018_25_6} (see Section 8). 

Define the value
$$\xi_n(\Omega):=\min \{ \xi(\Omega;S): \,
S  \mbox{ is an $n$-dimensional simplex,} \,
S\subset \Omega, \, \vo(S)\ne 0\}.$$
In other words, $\xi_n(\Omega)$ is the minimal absorption index of a convex body $\Omega$ by an inner nondegenerate simplex.  The results of the present survey are connected mainly with the case $\Omega=Q_n$. 

Let
$\xi_n:=\xi_n(Q_n).$
If $S\subset Q_n$, then $\xi(S)\geq n$, see
\eqref{xi_S_alpha_S_geq_n}.
This estimate occurs to be exact in order of   $n$. 
Really, consider for  $n>2$  the simplex with the zero-vertex and other $n$ vertices coinciding with the  vertices of $Q_n$ adjacent to $(1,\ldots,1)$.
 For this simplex, 
$\xi(S)=\frac{n^2-3}{n-1}.$ 
Hence, for $n>2$,  
\begin{equation}\label{xi_n_leq_frac_general}
\xi_n\leq \frac{n^2-3}{n-1}.
\end{equation}
If $n>1$, the right-hand part of  \eqref{xi_n_leq_frac_general} is  is strictly smaller than  $n+1$. The
inequality $\xi_n<n+1$
holds true also for $n=1,2$. If a simplex$S\subset Q_n$
has the maximum volume in $Q_n$, then  $\xi(S)\leq n+2$. For details, see  \cite{nevskii_monograph}.

  Thus, always $n\leq \xi_n< n+1$, i.\,e., $\xi_n-n\in [0,1)$.
However,  the exact values  of  
$\xi_n$ still are known only for 
$n=2$, $n=5$, $n=9$ and also for infinite set of $n$ when there exists an Hadamard matrix of order $n+1$. In all
these cases except 
$n=2$, holds $\xi_n=n$  while $\xi_2=\frac{3\sqrt{5}}{5}+1=2.34\ldots$
Dimension $n=2$ still remains a unique from even $n$ 
when  $\xi_n$ was calculated exactly.

The accurate value of
$\xi_n(B_n)$  is given in \cite{nevskii_mais_2018_25_6}.
If $S\subset B_n$, then $\xi(B_n;S)\geq n$, with an equality only for a regular simplex inscribed into
the boundary sphere. \linebreak (For brevity, further we will say that such a simplex is inscribed into the ball.)
Consequently,   $\xi_n(B_n)=n$. 



Let $x^{(j)}\in\Omega$ be the vertices of an  $n$-dimensional nondegenerate simplex  $S$.  
We say that an interpolation projector 
 $P:C(\Omega)\to \Pi_1({\mathbb R}^n)$ corresponds to
$S$ if its interpolation nodes coincides with the points $x^{(j)}$.
This projector is given by the equalities
$Pf\left(x^{(j)}\right)=
f_j:=
f\left(x^{(j)}\right).$
The following analogue of the Lagrange interpolation formula takes place:
\begin{equation}\label{interp_Lagrange_formula}
Pf(x)=\sum\limits_{j=1}^{n+1}
f\left(x^{(j)}\right)\lambda_j(x). 
\end{equation}
Let us denote by $\|P\|_\Omega$ the norm of  $P$ as an operator from $C(\Omega)$
to $C(\Omega)$. It follows from (\ref{interp_Lagrange_formula}) that


$$\|P\|_\Omega=\sup_{\|f\|_{C(\Omega)}=1} \|Pf\|_{C(\Omega)}=
\sup_{-1\leq f_j\leq 1} \max_{x\in\Omega} \sum_{j=1}^{n+1} f_j\lambda_j(x).$$
The expression $\sum f_j\lambda_j(x)$ is linear in $x$ and
$f_1,\ldots,f_{n+1},$ hence,
\begin{equation}\label{norm_P_intro_cepochka}
\|P\|_\Omega=\max_{f_j=\pm 1} \max_{x\in\Omega} \sum_{j=1}^{n+1}
f_j\lambda_j(x)
=\max_{x\in\Omega}\max_{f_j=\pm 1} \sum_{j=1}^{n+1}
f_j\lambda_j(x)
=\max_{x\in\Omega}\sum_{j=1}^{n+1}
|\lambda_j(x)|.
\end{equation}
If $\Omega$ is a convex polytope in ${\mathbb R}^n$ (e.\,g., $\Omega=Q_n)$,  
then we have also a simpler equality 
$$\|P\|_\Omega=
\max_{x\in\ver(\Omega)}\sum_{j=1}^{n+1}
|\lambda_j(x)|.$$

By $\theta_n(\Omega)$ define the minimal value of  $\|P\|_\Omega$
under the condition \linebreak $x^{(j)}\in \Omega$. Put $\theta_n:=\theta_n(Q_n)$. 

We say that a point$x\in\Omega$ is 
{\it an $1$-point of $\Omega$ with respect to $S,$}
if the projector
$P:C(\Omega)\to\Pi_1\left({\mathbb R}^n\right)$ with nodes in the vertices of $S$ satisfies an equality
$\|P\|=\sum
|\lambda_j(x)|$
and among numbers $\lambda_j(x)$ there is a unique negative.
It was proved in  \cite{nevskii_mais_2008_15_3} that for any
$P$ and the corresponding $S$ 
\begin{equation}\label{nev_ksi_P_ineq}
\frac{n+1}{2n}\Bigl( \|P\|_\Omega-1\Bigr)+1\leq 
\xi(\Omega;S) \leq
\frac{n+1}{2}\Bigl( \|P\|_\Omega-1\Bigr)+1.
\end{equation}
If   in $\Omega$ there exists an $1$-point   with respect to  $S,$ then the right-hand inequality becomes an
equality. The latter property in equivalent form was proved in \cite{nevskii_mais_2008_15_3}; 
the~notion {\it an 1-vertex of the cube}
was introduced later in \cite{nevskii_mais_2009_16_1}.

Suppose $\Omega=Q_n$ or $\Omega=B_n$. If $S\subset \Omega$, then $\xi(\Omega;S)\geq n.$ Applying
a right-hand inequality from \eqref{nev_ksi_P_ineq}, we see that for any projector having nodes in
 $\Omega$ 
$$\|P\|_{\Omega}\geq 3 - \frac{4}{n+1}.$$
This implies, for all $n$,
$$\theta_n\geq 3-\frac{4}{n+1}, \qquad \theta_n(B_n)\geq 3-\frac{4}{n+1}.$$ 

More advanced estimates of numbers  $\theta_n$ and $\theta_n(B_n)$ were obtained with the use of
the Legendre polynomials. {\it The standardized Legendre polynomial of degree  $n$} \linebreak is a function
$$\chi_n(t):=\frac{1}{2^nn!}\left[ (t^2-1)^n \right] ^{(n)}$$
(Rodrigues formula). For  properties of $\chi_n$ see, e.\,g., \cite{sege_1962}.
The Legendre polynomials are orthogonal upon the segment
$[-1,1]$ with the weight
$w(t)=1.$ 
As is known, $\chi_n(1)=1$;  if $n\geq 1$, then $\chi_n(t)$ 
increases for~$t\geq 1$.
By
$\chi_n^{-1}$, we denote a function inverse to $\chi_n$ 
on the  half-axis $[1,+\infty)$.

The appearance of the Legendre polynomials in our questions is connected with the following their property.
For $\gamma\geq 1$, let us consider the set 
$$E_{n,\gamma}:= \Bigl\{ x\in{\mathbb R}^n :
\sum_{j=1}^n\left |x_j\right| +
\Bigl|1- \sum_{j=1}^n x_j\Bigr|  \leq \gamma \Bigr\}.$$
In 2003, M.\,Nevskii \cite{nevskii_mais_2003_10_1} established the equality
\begin{equation}\label{mes_E_eq}
\mes(E_{n,\gamma}) 
= \frac{\chi_n(\gamma)}{n!} 
\end{equation}
(the proof is given also in \cite{nevskii_monograph}).
Applying this equality, the first author obtained  a~relation
\begin{equation}\label{theta_n_nu_h_ineq}
\theta_n
\geq
\chi_n^{-1} \left(\frac{1}{\nu_n}\right)=
\chi_n^{-1} \left(\frac{n!}{h_n}\right).
\end{equation}
Here $\nu_n$ is the maximum volume of a simplex contained in $Q_n$, 
$h_n$ is the maximum value of a determinant of order  $n$ consisting of $0$'s and $1$'s.
From 
(\ref{theta_n_nu_h_ineq}), some more visible inequalities were obtained.
This led to a relation $\theta_n\asymp \sqrt{n}.$  If $P:C(Q_n)\to \Pi_1({\mathbb R}^n)$ is an interpolation
projector having the nodes in vertices of a maximum volume simplex in  $Q_n,$ 
then $\|P\|_{Q_n}\asymp \theta_n$ (with constants not depenging on $n$).

Let us mark here some inequalities (links and proofs  are 
available
 in  \cite{nevskii_monograph}). 
For $n\not=2$, 
$$\frac{1}{e}\sqrt{n-1}<\theta_n \leq \min \left( \frac{n+1}{2},
\frac{4\sqrt{e}}{3} \sqrt{n+2}+1 \right).
$$
The left-hand inequality is valid for an arbitrary $n$. Also 
$$\frac{1}{4}\sqrt{n}<\theta_n<3\sqrt{n}, \quad n\in{\mathbb N}. $$
For a projector corresponding to a maximum volume simplex, we have 
$$\|P\|_{Q_n}\leq \min\left( n+1, \frac{4\sqrt{e}}{3} \sqrt{n+2}+1 \right),$$
If additionally  $n+1$ is an Hadamard number, then  $\|P\|_{Q_n}\leq \sqrt{n+1}$.

An analogue of \eqref{theta_n_nu_h_ineq} holds true  also for interpolation by linear functions on a~Euclidean ball. It is proved in  \cite{nevskii_mais_2019_26_3} that for each $n$
\begin{equation}\label{theta_n_chi_n_kappa_sigma_ineq_intro}
\theta_n(B_n)
\geq
\chi_n^{-1} 
\left(\frac{\varkappa_n}{\sigma_n}\right).
\end{equation}
Here  $\varkappa_n:=\vo(B_n)$ and 
$\sigma_n$ denotes the volume of an $n$-dimensional  regular simplex inscribed into   $B_n$. 
Utilizing \eqref{theta_n_chi_n_kappa_sigma_ineq_intro}, we obtain a relation $\theta_n(B_n)\asymp \sqrt{n}.$
 If $P:C(B_n)\to \Pi_1({\mathbb R}^n)$  is an interpolation projector with the nodes in vertices of a regular simplex inscribed into $B_n$, then
$\sqrt{n}\leq \|P\|_{B_n}\leq \sqrt{n+1}$.
 In other words, for such a projector 
 $\|P\|_{B_n}\asymp \theta_n(B_n)$  (see~\cite{nevskii_mais_2019_26_3}). 
 Moreover, at least for $1\leq n\leq 4$ we have $\|P\|_{B_n}= \theta_n(B_n)$ 
 (this is proved in \cite{nev_ukh_mais_2019_26_2} and \cite{nevskii_mais_2021_28_2}
 by different methods).

Nowadays, the exact values of
$\theta_n$ are known only for $n=1, 2, 3,$ and  $n=7$ (see~\cite{nevskii_monograph}). Namely,
$$\theta_1=1, \quad \theta_2=\frac{2\sqrt{5}}{5}+1=1.89\ldots, \quad \theta_3=2, \quad
\theta_7=\frac{5}{2}.$$
For all these dimensions, 
$$\xi_n=
\frac{n+1}{2}\Bigl( \theta_n-1\Bigr)+1.$$
Moreover, if $n=1,3, 7$, then
$$\theta_n= 3-\frac{4}{n+1}.$$

The exact values of $\theta_n(B_n)$ are known for $1\leq n\leq 4$\,:
$$\theta_1(B_1)=1, \quad \theta_2(B_2)=\frac{5}{3}, \quad \theta_3(B_3)=2, \quad \theta_4(B_4)=\frac{11}{5}.$$
In all these cases
$$\theta_n(B_n)= 3-\frac{4}{n+1}.$$ 
If $n\geq 5$, then  
$\theta_n(B_n)>3-\frac{4}{n+1}$. These results are obtained in  \cite{nev_ukh_mais_2019_26_2}.

In the present survey we will also discuss some estimates of numbers
$\xi_n^\prime$ и $\theta_n^\prime$ which are the analogues of 
$\xi_n$ and $\theta_n$ under the condition when 
when vertices of simplices coincide with vertices of the cube. By definition,
$$\xi_n^\prime:=\min\{ \xi(S) : \ver(S)\subset \ver(Q_n)\}, \quad
 \theta_n^\prime:=\min\{ \|P\|_{Q_n} : \ver(S)\subset \ver(Q_n)\}.$$
 In  the second equality 
  $S$ is  a  simplex with vertices in the nodes of
 a projector $P:C(Q_n)\to \Pi_1({\mathbb R}^n)$. Obviously, 
$$\xi_n^\prime\geq  \xi_n, \quad \theta_n^\prime\geq  \theta_n,$$
 This leads to the corresponding inequalities for upper and lower bounds of these values. 
 The numbers $\theta_n^\prime$ were introduced  in  paper \cite{nevskii_mais_2003_10_1} which is
 one of the very first M.\,Nevskii's works in this field. The estimates of $\xi_n^\prime$  and $\theta_n^\prime$ 
 are collected  iп \cite{nev_ukh_saratov_2018}; we will describe these results  in Section 4.
 
 For calculation and graph drawing,  the authors often
utilized   the Wolfram Ma\-the\-matica system (see, e.\,g.,  \cite{wolfram_2016}, \cite{klimov_uhalov_wolfram_2014}, \cite{ukhalov_practicum_2020}).


\section{The Case when  $n+1$ is  an Hadamard Number}\label{nev_ukh_sec_2}

%

{\it An Hadamard matrix of order   $m$} is a nondegenerate 
$(m\times m)$-matrix ${\bf H}_m$ with entries $\pm 1$ which satisfies
$${\bf H}_m^{-1}=
\frac{1}{m}{\bf H}_m^T.$$
Some information on Hadamard matrices is given in M.\,Hall's book~\cite{hall_1970}. 
If
${\bf H}_m$ exists, then $m=1,$ $m=2$ or $m$ is a multiple of $4.$ 
For an infinite set of numbers having form $m=4k$, 
including $m=2^l,$ the existence of 
${\bf H}_m$ is  long ago established. 
By~2008, the~smallest $m$, for which it is not known whether there exists an Hadamard matrix of order $m$,
was equal to $668.$ 
We will call $m$ {\it an Hadamard number} if there exists an~Hadamard matrix of order~$m$.

%

When and only when $n+1$ is  an Hadamard number, there exists a regular simplex $S$ which is 
inscribed into $Q_n$ in such a way that the vertices of $S$ coincide with vertices of the cube (see 
\cite[Theorem 4.5]{hudelson_1996}). 
In   \cite{nevskii_mais_2011_2} 
and  \cite[\S\,3.2]{nevskii_monograph},  it was proved  by various methods that 
$\xi(S)=n,$ so $\xi_n=n$. 
Here we give another proof given in \cite{nev_ukh_beitrage_2018}
and directly related to Hadamard matrices.


\smallskip
{\bf Theorem  2.1.}
{\it
Let $n+1$ be an Hadamard number. Then there exists a regular simplex $S$
with the following properties:

\smallskip
\noindent $1)$ $\ver(S) \subset \ver(Q_n)$; 

\smallskip
\noindent $2)$ $\xi_n=\xi(S)=\alpha(S)=n$;

\smallskip
\noindent $3)$ $d_1(S)=\ldots =d_n(S)=1$;

\smallskip
\noindent $4)$   the simplex $nS$ is circumscribed around $Q_n$ in such a way that
each $(n-1)$-dimensional face of $nS$ contains the only one vertex of $Q_n$.
}

\smallskip 
\smallskip
{\it Proof.}
Taking into consideration similarity, we can prove the statement for the cube $Q^\prime_n=[-1,1]^n$.
Since $n+1$ is an Hadamard number, there exists  {\it the~normalized Hadamard matrix} of order $n+1$, 
i.\,e., the Hadamard matrix with the  first row and the first column consisting of $1$'s (see
\cite[Chapter\,14]{hall_1970}).
Let us write the rows of this matrix in the inverse order:
$${\bf H} =
\left( \begin{array}{ccccc}
1&1&1&\ldots&1\\
\ldots&\ldots&\ldots&\ldots&1\\
\ldots&\ldots&\ldots&\ldots&1\\
\ldots&\ldots&\ldots&\ldots&\ldots\\
\ldots&\ldots&\ldots&\ldots&1\\
\end{array}
\right).$$
The obtained matrix {\bf H} also is an Hadamard matrix of order $n+1$. 
Consider the simplex $S^\prime$ with vertices 
formed by the first $n$ numbers in rows of $\bf H$.

All the vertices of $S^\prime$ are also vertices of $Q^\prime_n$, therefore,
the simplex is inscribed into the cube. Let us show that $S^\prime$ is a regular simplex
and the length of its edge is equal to $\sqrt{2(n+1)}$. 
Let  $a$, $b$ be two different rows of ${\bf H}$. All elements of 
$\bf H$ are $\pm 1$, hence, we have 
$\|a\|^2=\|b\|^2=n+1$. 
Rows of an Hadamard matrix are mutually orthogonal, therefore,
$$
\|a-b\|^2=(a-b,a-b)=\|a\|^2+\|b\|^2-2(a,b)=2(n+1).
$$
Denote by $u$ and $w$ the vertices of  $S^\prime$ obtained from  $a$ and $b$
respectively by removing the last component. 
This last component is equal to 1. It follows that $n$-dimensional length of vector $u-w$ is 
equal to $(n+1)$-dimensional length  of vector $a-b$, i.\,e., $\|u-w\|=\sqrt{2(n+1)}$. 

Denote by $\lambda_j$ the basic Lagrange polynomials of $S^\prime$. 
Since ${\bf H}^{-1}=
\frac{1}{n+1}{\bf H}^T$, 
the coefficients of 
$(n+1)\lambda_j$ are situated in rows of  ${\bf H}$. 
All the constant terms of these polynomials stand in the last column of ${\bf H}$. 
Consequently, they are equal to $1$. This means that the constant terms of polynomials
$-(n+1)\lambda_j$ are equal to $-1$.
Hence, for any $j=1,\ldots,n+1$, 
$$(n+1)\max_{x\in\ver(Q^\prime_n)}(-\lambda_j(x))=n-1.$$
The coefficients of polynomials $-(n+1)\lambda_j$ are equal to $\pm 1$, 
therefore, for any $j$, the vertex  $v$ of $Q^\prime_n$, such that $(n+1)(-\lambda_j(v))=n-1$,
is unique.
Indeed, $v=(v_1,\ldots,v_n)$ is defined by equalities $v_i=- \operatorname{sign} l_{ij}$, where  $l_{ij}$ are the  coefficients of $\lambda_j$. 

Now let us find $\xi(Q^\prime_n; S^\prime)$ using 
a formula analogous to (\ref{xi_S_formula}): 
$$\xi(Q^\prime_n;S^\prime)=(n+1)\max_{1\leq j\leq n+1}
\max_{x\in \ver(Q^\prime_n)}(-\lambda_j(x))+1=n-1+1=n.
$$
Consider the similarity transformation which maps $Q^\prime_n$ into $Q_n$.
This transformation also maps $S^\prime$ into a simplex inscribed into $Q_n$. 
Denote by  $S$ the image of $S^\prime$.
Obviously, $\xi(S)=\xi(S^\prime;Q^\prime_n)=n$.
It follows that, if $n+1$ is an Hadamard number, then $\xi_n\leq n$. 
As we know, for any $n$, $\xi_n\geq n$.
Hence, $\xi_n=\xi(S)=n.$ 
Since the~coefficient of similarity for a  mapping $Q^\prime_n$
to $Q_n$ equals $\frac{1}{2}$, the length of any edge of $S$ equals $\sqrt{   \frac{n+1}{2}  }$.

The condition
$$
\max\limits_{x\in \ver(Q^\prime_n)} \left(-\lambda_1(x)\right)=
\ldots=
\max\limits_{x\in \ver(Q^\prime_n)} \left(-\lambda_{n+1}(x)\right)=\frac{n-1}{n+1}
$$
and \eqref{xi_S_S_circ_condition_around_Q_n} mean that the simplex $nS^\prime$ 
is circumscribed around the cube $Q^\prime_n$. 
From the above it also follows that
each $(n-1)$-face of $nS^\prime$ 
contains the only one vertex of $Q^\prime_n$.
Hence, the simplex $nS$ is circumscribed around the cube $Q_n$ in the same way.
Since $nS$ is circumscribed around $Q_n$, we have $\alpha(S)=\xi(S)=n$
and $d_i(S)=1$. These equalities could be also obtained from inequalities
$n\leq \alpha(S)\leq \xi(S)$ and the~equality $\xi(S)=n$. 

The equalities  $d_i(S)=1$ and $\alpha(S)=n$ can be also derived in another way.
Note that a regular simplex  inscribed into $Q_n$ has the maximum volume 
among all  the simplices contained  in this cube   (see  \cite[Theorem~2.4]{hudelson_1996}).
As we marked, each axial diameter of any maximum volume simplex in $Q_n$  is equal to $1$.

The theorem is proved.
\hfill$\Box$


\section{Estimates for the Minimal Absorption Index \\ of   
a Cube by a  Simplex
}
\label{nev_ukh_sec_3}


The above estimates
\begin{equation}\label{xi_n_leq_rough}
n\leq \xi_n< n+1, \quad n\in{\mathbb N}
\end{equation}
and
$$
\xi_n\leq \frac{n^2-3}{n-1}, \quad n>2 
$$
 were  recently specified for some dimensions $n$.
Here we give all the exact values and the best estimates of $\xi_n$ known for today. For $n=1,\ldots,10$,
see Table ~\ref{tab:nev_uhl_ksi_est_1_10}.

\begin{table}[!htbp]
\begin{center}
\caption{ \label{tab:nev_uhl_ksi_est_1_10}
Estimates of $\xi_n$\  $(1\leq n\leq 10)$ }
\bgroup
\def\arraystretch{2.1}%
\begin{tabular}{|c|c|c|c|}
\hline
$n$ & $\xi_n$ &$n$ & $\xi_n$\\
\hline
 $1$ & $1$ &  $6$ & $6\leq \xi_6< 6.0166$\\
 \hline
 $2$  & $\frac{3\sqrt{5}}{5}+1=2.3416\ldots$  &      $7$   & $7$   \\
\hline
 $3$  & $3$  &      $8$   &  $8\leq \xi_8 < 8.1355$ \\
\hline
 $4$  & $4\leq \xi_4\leq \frac{19+5 \sqrt{13}}{9}=4.1141\ldots$  & $9$  & $9$   \\
\hline
 $5$  & $5$   &      $10$   & $\qquad\  10\leq \xi_{10} < 10.2342\ \qquad$  \\
\hline
\end{tabular}
\egroup
\end{center}
\end{table}

Clearly, $\xi_1=1$. The exact values of $\xi_2$, $\xi_3$, and $\xi_7$ were  delivered  by M.\,Nevskii  (see~\cite{nevskii_monograph}).
Note that $\xi_1$, $\xi_3$, and $\xi_7$ can be also calculated as a corollary of Theorem 2.1, since $2, 4$,
and $8$ are Hadamard numbers.
As is noted, $n=2$ remains a unique even $n$ for which we  know the exact $\xi_n$ .

The given estimate for  $\xi_4$ was obtained in \cite{nev_ukh_mais_2016_23_5}. 
In this paper, also an assumption was made that 
$\xi_4=\frac{19+5 \sqrt{13}}{9}$.
The authors have not yet succeeded in proving or refuting this conjecture.

The value of $\xi_5$ and the clarified estimate for  $\xi_6$ were firstly given in \cite{nev_ukh_mais_2017_24_1}.
It~was shown that there exists a one-parameter family of simplices   $R=R(t)$ with a property $\xi(R)=5$. 
In the same way, in \cite{nev_ukh_beitrage_2018} it was established that  $\xi_9=9$.

For today, $n=5$ is the minimal $n$ for which an infinite set of extremal simplices was found.
Later such the families were built also for $n=7$ and $n=9$ (see~\cite{nev_ukh_beitrage_2018}).
For $n=5$,  in~\cite{nev_ukh_beitrage_2018} a two-parameter family of simplices $R=R(s,t)$ satisfying $\xi(R)=5$
was built.  In addition to being extremal, simplices from this family also have some other interesting properties (see Sections 6--8). 

The most precise upper bounds for $\xi_8$ and $\xi_{10}$ were obtained by
A.\,Ukhalov and D.\,Fedulov in    \cite{fedulov_ukhalov_2019}. 

The given estimates for $\xi_6$, $\xi_8$, and $\xi_{10}$ were discovered by numerical minimization
of $\xi(S)$ as a function of  vertices coordinates of a simplex.

\begin{table}[!htbp]
\begin{center}
\caption{\label{tab:nev_uhl_ksitab1} Estimates of $\xi_n$\  $(11\leq n \leq 58)$}
\bgroup
$
\def\arraystretch{2.4}
\begin{array}{||c|c|c||c|c|c||c|c|c|}
\hline
 {\bf n} &  {\bf \xi_n\leq}  &  & {\bf n} &   {\bf\xi_n\leq}  &  &   {\bf n} &   {\bf\xi_n\leq}  & \\
 \hline
 11 & 11 & \text{H} & 27 & 27  & \text{H} & 43 & 43  & \text{H} \\ \hline
 12 & \frac{184}{15} = 12.267 & \text{B} & 28 & \frac{781}{27} = 28.926 & \text{A} & 44 & \frac{1933}{43} = 44.953 & \text{A} \\ \hline
 13 & \frac{529}{39} = 13.564 & \text{B} & 29 & \frac{5963}{203} = 29.374 & \text{B} & 45 & \frac{2493}{55} = 45.327 & \text{B} \\ \hline
 14 & \frac{193}{13} = 14.846 & \text{A} & 30 & \frac{2981}{98}= 30.418 & \text{B} & 46 & \frac{2113}{45} = 46.956 & \text{A} \\ \hline
 15 & 15  & \text{H} & 31 & 31  & \text{H} & 47 & 47  & \text{H} \\ \hline
 16 & \frac{163}{10} = 16.300 & \text{B} & 32 & \frac{1600}{49} = 32.653 & \text{B} & 48 & \frac{2301}{47} = 48.957 & \text{A} \\ \hline
 17 & \frac{296}{17} = 17.412 & \text{B} & 33 & \frac{543}{16} = 33.938 & \text{A} & 49 & \frac{689}{14} = 49.214 & \text{B} \\ \hline
 18 & \frac{321}{17} = 18.882 & \text{AB} & 34 & \frac{5237}{152} = 34.454 & \text{B} & 50 & \frac{1162}{23} = 50.522 & \text{B} \\ \hline
 19 & 19  & \text{H} & 35 & 35  & \text{H} & 51 & 51  & \text{H} \\ \hline
 20 & \frac{596}{29} = 20.552 & \text{B} & 36 & \frac{875}{24} = 36.458 & \text{B} & 52 & \frac{2701}{51} = 52.961 & \text{A} \\ \hline
 21 & \frac{219}{10} = 21.900 & \text{A} & 37 & \frac{4139}{111} = 37.288 & \text{B} & 53 & \frac{36707}{689} =53.276 & \text{B} \\ \hline
 22 & \frac{481}{21} = 22.905 & \text{A} & 38 & \frac{2309}{60} =38.483 & \text{B} & 54 & \frac{2017}{37} = 54.514 & \text{B} \\ \hline
 23 & 23  & \text{H} & 39 & 39  & \text{H} & 55 & 55  & \text{H} \\ \hline
 24 & \frac{339}{14} = 24.214 & \text{B} & 40 & \frac{1808}{45} = 40.178 & \text{B} & 56 & \frac{394}{7} = 56.286 & \text{B} \\ \hline
 25 & \frac{379}{15} = 25.267 & \text{B} & 41 & \frac{8479}{205} = 41.361 & \text{B} & 57 & q = 57.618 & \text{B} \\ \hline
 26 & \frac{4853}{182} = 26.665 & \text{B} & 42 & \frac{40204}{945} = 42.544 & \text{B} & 58 & \frac{31455}{539} = 58.358 & \text{B} \\ \hline
\end{array}
$
\egroup
\end{center}
\end{table}

\begin{table}[!htbp]
\begin{center}
\caption{\label{tab:nev_uhl_ksitab2} Estimates of $\xi_n$\  $(59\leq n \leq 106)$}
\bgroup
$
\def\arraystretch{2.4}
\begin{array}{||c|c|c||c|c|c||c|c|c|}
\hline
 {\bf n} &  {\bf \xi_n\leq}  &  & {\bf n} &   {\bf\xi_n\leq}  &  &   {\bf n} &   {\bf\xi_n\leq}  & \\
 \hline
 59 & 59  & \text{H} & 75 & 75  & \text{H} & 91 & 91  & \text{H} \\ \hline
 60 & \frac{1985}{33} = 60.152 & \text{B} & 76 & \frac{5773}{75} = 76.973 & \text{A} & 92 & \frac{8461}{91} = 92.978 & \text{A} \\ \hline
 61 & \frac{11219}{183} = 61.306 & \text{B} & 77 & \frac{2963}{38} = 77.974 & \text{A} & 93 & \frac{4323}{46} = 93.978 & \text{A} \\ \hline
 62 & \frac{21299}{341} = 62.460 & \text{B} & 78 & \frac{4081}{52} = 78.481 & \text{B} & 94 & \frac{5857}{62} = 94.468 & \text{B} \\ \hline
 63 & 63  & \text{H} & 79 & 79  & \text{H} & 95 & 95  & \text{H} \\ \hline
 64 & \frac{4093}{63} = 64.968 & \text{A} & 80 & \frac{3937}{49} = 80.347 & \text{B} & 96 & \frac{9213}{95} = 96.979 & \text{A} \\ \hline
 65 & \frac{33949}{520} = 65.287 & \text{B} & 81 & \frac{3653}{45} = 81.178 & \text{B} & 97 & \frac{18863}{194} = 97.232 & \text{B} \\ \hline
 66 & \frac{5918}{89} = 66.494 & \text{B} & 82 & \frac{8990}{109} = 82.477 & \text{B} & 98 & \frac{4234}{43}= 98.465 & \text{B} \\ \hline
 67 & 67  & \text{H} & 83 & 83  & \text{H} & 99 & 99  & \text{H} \\ \hline
 68 & \frac{4621}{67} = 68.970 & \text{A} & 84 & \frac{42406}{501} = 84.643 & \text{B} & 100 & \frac{714817}{7110} = 100.537 & \text{B} \\ \hline
 69 & \frac{2379}{34} = 69.971 & \text{A} & 85 & \frac{152113}{1785} = 85.217 & \text{B} & 101 & \frac{51097}{505} = 101.182 & \text{B} \\ \hline
 70 & \frac{3313}{47} = 70.489 & \text{B} & 86 & \frac{1643}{19} = 86.474 & \text{B} & 102 & \frac{19555}{191} = 102.382 & \text{B} \\ \hline
 71 & 71  & \text{H} & 87 & 87 & \text{H} & 103 & 103  & \text{H} \\ \hline
 72 & \frac{5181}{71} = 72.972 & \text{A} & 88 & \frac{7741}{87} = 88.977 & A^* & 104 & \frac{10813}{103} =104.981 & A^* \\ \hline
 73 & \frac{48128}{657} = 73.254 & \text{B} & 89 & \frac{3959}{44} = 89.977 & A^* & 105 & \frac{5511}{52} = 105.981 & \text{A} \\ \hline
 74 & \frac{2458}{33}= 74.485 & \text{B} & 90 & \frac{1538}{17} = 90.471 & \text{B} & 106 & \frac{7021}{66} = 106.379 & \text{B} \\ \hline
\end{array}
$
\egroup
\end{center}
\end{table}

\begin{table}[!htbp]
\begin{center}
\caption{\label{tab:nev_uhl_ksitab3} Estimates of $\xi_n$\ $(107\leq n \leq 118)$}
\bgroup
$
\def\arraystretch{2.4}
\begin{array}{||c|c|c||c|c|c|}
\hline
 {\bf n} &  {\bf \xi_n\leq}  &  & {\bf n} &   {\bf\xi_n\leq}  &   \\
 \hline
 107 & 107  & \text{H} & 114 & \frac{24247}{212} = 114.373 & \text{B} \\ \hline
 109 & \frac{107131}{981}= 109.206 & \text{B} & 115 & 115  & \text{H} \\ \hline
 110 & \frac{22627}{205} = 110.376 & \text{B} & 116 & \frac{13453}{115} = 116.983 & A^* \\ \hline
 111 & 111  & \text{H} & 117 & \frac{132583}{1131} = 117.226 & \text{B} \\ \hline
 112 & \frac{6727}{60} = 112.117 & \text{B} & 118 & \frac{8641}{73} = 118.370 & \text{B} \\ \hline
 113 & \frac{25591}{226} = 113.235 & \text{B} & \text{} & \text{} & \text{}  \\ \hline
\end{array}
$
\egroup
\end{center}
\end{table} 

Tables \ref{tab:nev_uhl_ksitab1}, \ref{tab:nev_uhl_ksitab2}, and \ref{tab:nev_uhl_ksitab3}
contain the most precise upper bounds of $\xi_n$ for $11\leq n\leq 118$. 
These estimates are obtained in \cite{nev_ukh_mais_2018_25_1}. In order to provide the estimates, we used
the following approach. M.\,Nevskii \cite{nevskii_monograph} showed that if $S\subset Q_n$ is a maximum
volume simplex,then $\xi(S)\asymp \xi_n$. 
Calculation of $\xi(S)$ for maximum volume simplices often
makes it possible to  refine the upper bound in
the double inequality 
\begin{equation}\label{xi_n_double_ineq}
n\leq \xi_n\leq \frac{n^2-3}{n-1}
\end{equation}
(see \eqref{xi_n_leq_frac_general} and \eqref{xi_n_leq_rough}). For constructing
maximum volume simplices  in $Q_n$, we utilized some data from website
http://www.indiana.edu/\char`\~maxdet/  \quad
related to maximum  $(-1,1)$-determinants.
This site contained maximum $(-1,1)$-determinants up to the order $119$. 
Among the listed determinants
there were also those whose maximality has not been established. For some orders, the determinants
were not given at all.

Along with  upper bounds of $\xi_n$, the tables also provide  information about
 methods of obtaining the estimates. We use the following designations:


\begin{itemize}
\item A --- The estimate is obtained from inequality \eqref{xi_n_double_ineq}. 
The using of a maximum volume simplex did not improve this estimate. 

\item B ---  We give the smallest value $\xi(S)$ obtained from considering one 
or more simplices of maximum volume in $Q_n.$ 
This value is less than the right-hand side of \eqref{xi_n_double_ineq}. 

\item AB --- The value of $\xi(S)$ obtained for simplices of maximum volume 
coincides with the right-hand side of \eqref{xi_n_double_ineq}.

\item H ---  $n+1$ is an Hadamard number. We give the exact value $\xi_n = n.$  

\item $\text{A}^*$ --- We use inequality \eqref{xi_n_double_ineq},
since for this $n$ no maximum determinant is given on the site.
\end{itemize}

For $n=57,$ we use denotation
$q=\frac{34018994107517188105}{590424166322794597}.$

%
%
%
%
%
%
%
%
%
%
%


\section{Estimates for the Minimal Norm of a Projector\\
in Linear Interpolation on a Cube in ${\mathbb R}^n$}\label{nev_ukh_sec_4}

Evidently,  $\theta_1=1$. The exact values of $\theta_2$, $\theta_3$, and $\theta_7$  were obtained 
 by M.\,Nevskii (see~\cite{nevskii_monograph}). For today, $n=1,$ $n=2$, $n=3$, and $n=7$
 are the only cases when we know $\theta_n$ exactly.
  For other $n$,  we know common estimates  (see Section 1).
 In some cases, these estimates were clarified with the use of computer. Here we mostly give the material 
 from \cite{nev_ukh_mais_2018_25_3}. 
 
 Table~\ref{tab:nev_uhl_theta_upper} contains the most precise upper estimates of 
 $\theta_n$ for $1\leq n \leq 27$.
  If we know the exact value, we write this value. For providing estimates, the following approach
 described in \cite{nevskii_monograph} may be used.
 The upper bound of  $\theta_n$  appears from consideration of an projector having the nodes in vertices
 of a maximum volume simplex in $Q_n$. For the first time, this method was applied by M.\,Nevskii and I.\,Hlestkova in  \cite{nev_hlestkova_2008}.

A significant part of the given estimates   was obtained using actually this approach.
For constructing maximum volume simplices for a given $n$, we utilized maximal $(-1,1)$-determinants of order $n+1$ from the site \newline
http://www.indiana.edu/\char`\~maxdet/ .
Hadamard matrices were imported from website 
http://www.maths.gla.ac.uk/\char`\~es/hadamard/hadamard.php .

\begin{table}[!htbp]
\begin{center}
\caption{\label{tab:nev_uhl_theta_upper} Upper estimates of $\theta_n$\ $(1\leq n \leq 27)$}
\bgroup
$
\def\arraystretch{1.5}
\begin{array}{|r|c|c|c|c|}
  \hline
 {\bf n} & {\bf min\, \|P\|} &  \text{{\bf N}} &   {\bf \theta_n\leq}  & \\
  \hline
 1 & 1 & 1 & \theta_1=1 & \text{N} \\
  \hline
 2 & 3 & 1  & \theta_2=\frac{2\sqrt{5}}{5}+1=1.8944\ldots & \text{N} \\
  \hline
 3 & 2 & 1  & \theta_3=2 &\text{N} \\
  \hline
 4 & \frac{7}{3}  =2.3333\ldots & 1  & \frac{3  (4+\sqrt{2})}{7}  =2.3203\ldots & \text{K} \\
  \hline
 5 & \frac{13}{5}= 2.6 & 1 & 2.448804 & \text{NU} \\
  \hline
 6 & 3  & 1  & 2.60004\ldots  & \text{L} \\
  \hline
 7 & \frac{5}{2} = 2.5 &  1  & \theta_7=\frac{5}{2} = 2.5 & \text{N} \\
  \hline
 8 & \frac{22}{7} = 3.1428\ldots & 1 & \frac{22}{7} = 3.1428\ldots & \text{M} \\
  \hline
 9 & 3 & 1 &  3 & \text{M} \\
  \hline
 10 & \frac{19}{5} = 3.8 & 3 &\frac{19}{5} = 3.5186\ldots &  \text{L} \\
  \hline
 11 & 3 & 1 & 3 & \text{H} \\
  \hline
 12 & \frac{17}{5} = 3.4 & 1 & \frac{17}{5} = 3.4 & \text{M} \\
  \hline
 13 & \frac{49}{13} = 3.7692\ldots & 1 & \frac{49}{13} = 3.7692\ldots & \text{M} \\
  \hline
 14 & \frac{21}{5} = 4.2 & 1 & \frac{21}{5} = 4.0156\ldots & \text{L} \\
  \hline
 15 & \frac{7}{2} = 3.5 & 5 & \frac{7}{2} = 3.5 & \text{H} \\
  \hline
 16 & \frac{21}{5} = 4.2 & 3 & \frac{21}{5} = 4.2 & \text{M} \\ 
 \hline
 17 & \frac{139}{34} = 4.0882\ldots & 3 & \frac{139}{34} = 4.0882\ldots &  \text{M} \\
  \hline
 18 & \frac{95}{17} = 5.5882\ldots & 3 & 5.14006\ldots & \text{L} \\
  \hline
 19 & 4 & 3 & 4 & \text{H} \\
  \hline
 20 & \frac{137}{29} = 4.7241\ldots & 7 &  4.68879\ldots & \text{L} \\
  \hline
 21 & \frac{251}{50} = 5.02 & 1 &  \frac{251}{50} = 5.02 & \text{M?} \\
  \hline
 22 & \frac{1817}{335} =5.4238\ldots & 1 & \frac{1817}{335} =5.4238\ldots & \text{M?} \\
  \hline
 23 & \frac{9}{2} = 4.5 & 60 & \frac{9}{2} = 4.5 & \text{O} \\
  \hline
 24 & \frac{103}{21} = 4.9047\ldots & 2 & \frac{103}{21} = 4.9047\ldots & \text{M} \\
  \hline
 25 & 5 & 3 & 5 & \text{M} \\
  \hline
 26 & \frac{474}{91}= 5.2087\ldots & 1 & \frac{474}{91}= 5.2087\ldots & \text{M?} \\
  \hline
  
 27 & 5.0 & 487 & 5.0 & \text{U} \\
   \hline
\end{array}
$
\egroup
\end{center}
\end{table} 

The maximality of $(-1,1)$-determinants of order $1,\ldots, 21$
listed on website \linebreak  www.indiana.edu/$\sim$maxdet  is proved.
Moreover, the website contains completely all 
the maximum determinants of order $\leq 21$ (up to equivalence of matrices).
While getting the estimates, we used all these data.
The maximality of determinants of order 22, 23, and 27 has not been established. 
Such situations
are marked in the table with a question mark. It is unclear whether is exhaustive
the set of maximum determinants of order 25 and 26 given on the site.
To refine the upper estimate of $\theta_n,$ we examined
all the known maximum $(-1,1)$-determinants of order
$n + 1,$ since the corresponding projectors
can have different norms.

Calculations using the formula
\begin{equation}\label{nev_uhl_normPQn}
\|P\|_{Q_n} =
\max_{x\in\ver(Q_n)}\sum_{j=1}^{n+1}
|\lambda_j(x)|.
\end{equation}
reduce to looking through all $2^n$ vertices of $Q_n.$
For large $n,$ these calculations become rather time-consuming. 
This is the reason why we limited by investigation of cases $1\leq n\leq 27.$

In all the situations when $n + 1$ is an Hadamard number, we used the complete set of
Hadamard matrices of the appropriate order. In particular, to obtain the estimate
of $\theta_{23},$ we considered all the existing 60 Hadamard matrices of order 24.
This case was studied by A.\,Ukhalov and E.\,Ozerova \cite{kudr_ozerova_ukh_2017}.
The estimate of  $\theta_{27}$ was obtained by A.\,Ukhalov and E.\,Udovenko \cite{udovenko_2019}, \cite{ukhalov_udovenko_2020}. To get the result, it was necessary to consider 487 Hadamard matrices of order  28.


For $n=4,5, 6,10,14,18,$ and $20$, the estimates of  $\theta_n$  were made via minimization
of the function standing in the right-hand part  \eqref{nev_uhl_normPQn}. The links are given in the description
of the table.

The authors  conjecture that 
 $\theta_4=\frac{3  (4+\sqrt{2})}{7}  =2.320377\ldots$

In the first column of 
Table~\ref{tab:nev_uhl_theta_upper}, we put
dimension  $n.$ 
The second column contains the minimal projector norm 
obtained via consideration of maximum $(-1,1)$-determinants of order $n + 1.$
The third column shows number  $N$ of matrices of order $n + 1$ 
considered for obtaining the result written in the second column.
In the fourth column, the best nowaday upper bound of $\theta_n $ is given.
The last column of the table contains a comment on the way to get the top
boundary of $\theta_n $ indicated in the fourth column.
Here we use the following notation.

%
%
%
%
%
%
%

\begin{itemize}
\item H --- The case when $n+1$ is an Hadamard number.  There was used a 
regular simplex built from the corresponding Hadamard matrix.

\item K --- The estimate is obtained by I.\,Kudryavcev (see~\cite{kudr_ozerova_ukh_2017}).

\item L --- The estimate is obtained by A.\,Lutenkov in \cite{lutenkov_2018}. 

\item M ---  There was used a maximum volume simplex in $Q_n$ 
built from maximum $(-1,1)$-determinant of order $n + 1.$

\item M? --- There was used a simplex in $Q_n$ 
built from determinant of order $n + 1$ which is maximum from known nowaday.    
The global maximality of such a determinant yet is not proved.

\item N --- The exact value of $\theta_n$ is obtained by  M.\,Nevskii (see~\cite{nevskii_monograph}). 

\item NU ---The upper estimate is obtained by   M.\,Nevskii and A.\,Ukhalov  in~\cite{nev_ukh_mais_2017_24_1}.

\item O --- The upper estimate is obtained by E.\,Ozerova (see~\cite{kudr_ozerova_ukh_2017}).

\item U --- The upper estimate is obtained by E.\,Udovenko in \cite{udovenko_2019}.
\end{itemize}

\begin{table}[!htbp]
\begin{center}
\caption{\label{tab:nev_uhl_theta_lower1} Lower estimates of $\theta_n$\ $(1\leq n \leq 54)$}
\bgroup
$
\def\arraystretch{1.44}
\begin{array}{|c|c|c|c||c|c|c|c|}
\hline
{\bf n} & {\bf \chi_n^{-1}(\frac{1}{\nu_n}) } & {\bf 3-\frac{4}{n+1}} & \text{\bf Max} &
{\bf n} & {\bf \chi_n^{-1}(\frac{1}{\nu_n}) } & {\bf 3-\frac{4}{n+1}} & \text{\bf Max} \\
\hline
 1 & 1 & 1 & 1 & 28 & 2.2768 & 2.8621 & 2.8621 \\
 \hline
 2 & 1.291 & 1.6667 & 1.6667 & 29 & 2.3074 & 2.8667 & 2.8667 \\
 \hline
 3 & 1.2492 & 2 & 2 & 30 & 2.3487 & 2.871 & 2.871 \\
 \hline
 4 & 1.3478 & 2.2 & 2.2 & 31 & 2.3452 & 2.875 & 2.875 \\
 \hline
 5 & 1.4284 & 2.3333 & 2.3333 & 32 & 2.3955 & 2.8788 & 2.8788 \\
 \hline
 6 & 1.5018 & 2.4286 & 2.4286 & 33 & 2.4259 & 2.8824 & 2.8824 \\
 \hline
 7 & 1.4678 & 2.5 & 2.5 & 34 & 2.4642 & 2.8857 & 2.8857 \\
 \hline
 8 & 1.5626 & 2.5556 & 2.5556 & 35 & 2.4601 & 2.8889 & 2.8889 \\
 \hline
 9 & 1.6034 & 2.6 & 2.6 & 36 & 2.5019 & 2.8919 & 2.8919 \\
 \hline
 10 & 1.6699 & 2.6364 & 2.6364 & 37 & 2.5348 & 2.8947 & 2.8947 \\
 \hline
 11 & 1.6488 & 2.6667 & 2.6667 & 38 & 2.5722 & 2.8974 & 2.8974 \\
 \hline
 12 & 1.7086 & 2.6923 & 2.6923 & 39 & 2.5697 & 2.9 & 2.9 \\
 \hline
 13 & 1.7659 & 2.7143 & 2.7143 & 40 & 2.6056 & 2.9024 & 2.9024 \\
 \hline
 14 & 1.8211 & 2.7333 & 2.7333 & 41 & 2.641 & 2.9048 & 2.9048 \\
 \hline
 15 & 1.8108 & 2.75 & 2.75 & 42 & 2.6759 & 2.907 & 2.907 \\
 \hline
 16 & 1.8778 & 2.7647 & 2.7647 & 43 & 2.6747 & 2.9091 & 2.9091 \\
 \hline
 17 & 1.9156 & 2.7778 & 2.7778 & 44 & 2.7179 & 2.9111 & 2.9111 \\
 \hline
 18 & 1.965 & 2.7895 & 2.7895 & 45 & 2.743 & 2.913 & 2.913 \\
 \hline
 19 & 1.9587 & 2.8 & 2.8 & 46 & 2.7791 & 2.9149 & 2.9149 \\
 \hline
 20 & 2.0159 & 2.8095 & 2.8095 & 47 & 2.7756 & 2.9167 & 2.9167 \\
 \hline
 21 & 2.0588 & 2.8182 & 2.8182 & 48 & 2.8201 & 2.9184 & 2.9184 \\
 \hline
 22 & 2.1039 & 2.8261 & 2.8261 & 49 & 2.8413 & 2.92 & 2.92 \\
 \hline
 23 & 2.0958 & 2.8333 & 2.8333 & 50 & 2.8805 & 2.9216 & 2.9216 \\
 \hline
 24 & 2.1408 & 2.84 & 2.84 & 51 & 2.8729 & 2.9231 & 2.9231 \\
 \hline
 25 & 2.1847 & 2.8462 & 2.8462 & 52 & 2.9173 & 2.9245 & 2.9245 \\
 \hline
 26 & 2.2278 & 2.8519 & 2.8519 & 53 & 2.9362 & 2.9259 & 2.9362 \\
 \hline
 27 & 2.2242 & 2.8571 & 2.8571 & 54 & 2.9735 & 2.9273 & 2.9735 \\
 \hline
\end{array}
$
\egroup
\end{center}
\end{table} 

For $n=5$, E.\,Spasova (2020, bachelor's thesis under the guidance of M.\,Nevskii)  constructed a projector 
having the norm $\|P\|_{Q_5}=\frac{5}{2}=2.5.$ 
This value exceeds the smallest known  upper bound
$ 2.4488 \ldots $ for $ \theta_5 $ established in \cite{nev_ukh_mais_2017_24_1}.
Nevertheless, the nodes of  $P$ have a rather simple form ---  they are located at the vertices  or at midpoints of the edges of $Q_5$.
For example, the following  points are suitable:
$(0, 0, 0, 0, 0),  \left(0, 0, \frac{1}{2}, 1, 1 \right),$ $\left(0, 1, \frac{1}{2},0, 1\right),
(0, 1, 1,1, 0), \left(1,0,1,0, \frac{1}{2}\right),$ $ \left(1,1, 0, 1, \frac{1}{2}) \right).$ \linebreak
Back in 2006,   a projector with nodes in the vertices of $Q_5$ and the norm
$2.6$ was built (see, e.\,g., \cite{nevskii_mais_2006_13_2}), and for at least 10 years the estimate $\theta_5 \leq 2.6 $ was the best from the known.

Now let us pass to the lower bounds of $\theta_n$. Tables~\ref{tab:nev_uhl_theta_lower1}, \ref{tab:nev_uhl_theta_lower2}, and \ref{tab:nev_uhl_theta_lower3} present the largest known
lower bounds for these numbers. To obtain the estimates, we used the inequality
\begin{equation}\label{nev_theta_n_max_ineq}
\theta_n
\geq  \max \left[
3 - \frac{4}{n+1},
\chi_n^{-1} \left(\frac{1}{\nu_n}\right) \right].
\end{equation}
For each $n$, we give  the values  $ \chi_n^{-1}\left(\frac{1}{\nu_n}\right)$ and
$3-\frac{4}{n+1}$, as well as the maximum of them. This maximum is the most accurate lower bound
for $\theta_n$.
The tables are missing $ n = 104 $ and $ n = 116 $, since
 website www.indiana.edu/$\sim$maxdet does not show any maximum
determinants of orders 105 and 117.

\begin{figure}[!htbp]
 \centering
\includegraphics[width=14.7cm]{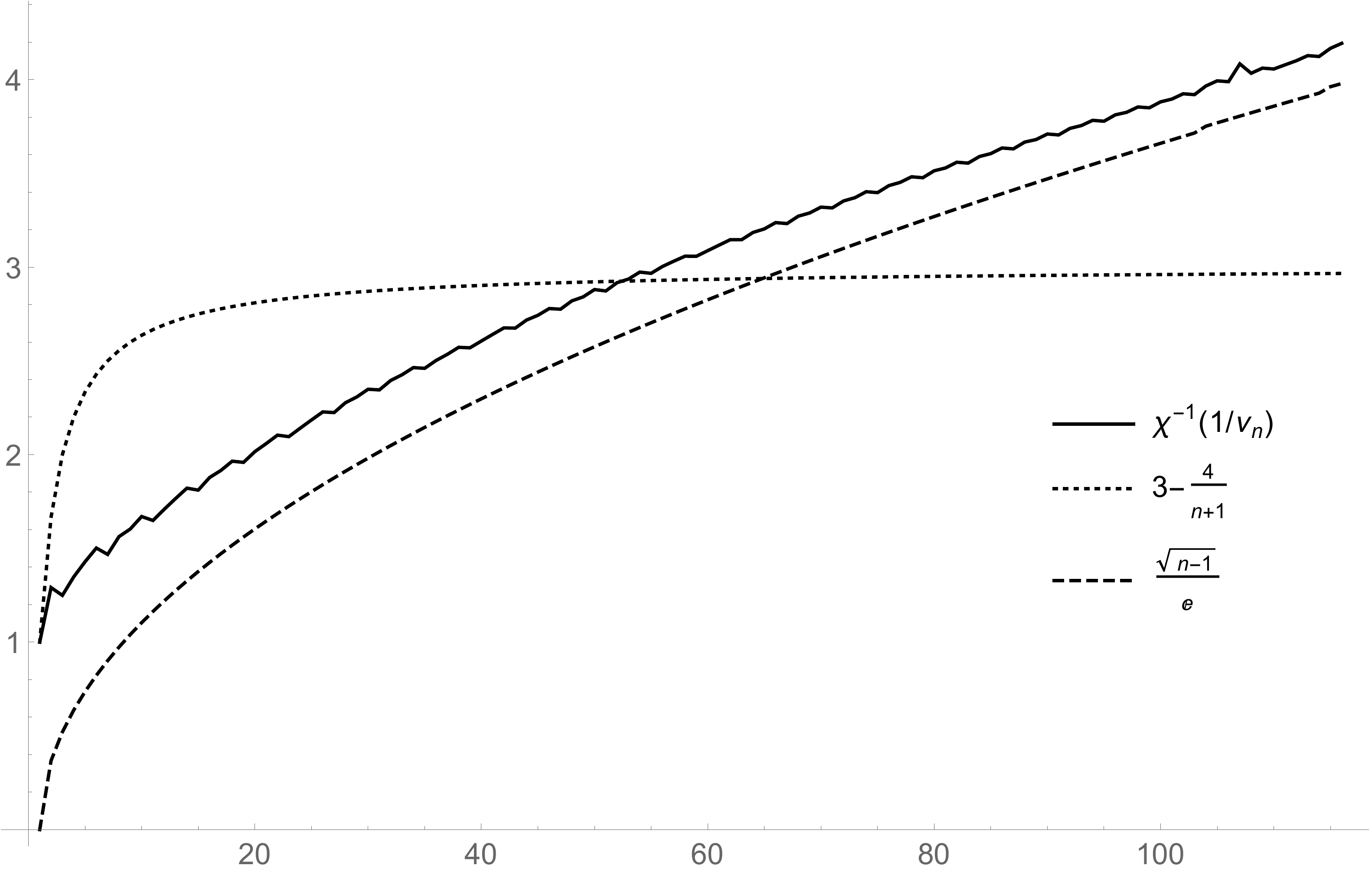}
\caption{Graphs of the functions $\chi_n^{-1}\left(\frac{1}{\nu_n}\right)$,  \,
$3-\frac{4}{n+1}$\,,  and \, $\frac{\sqrt{n-1}}{e}$}
\label{fig:nev_ukl_theta_lower}
\end{figure}

Graphs of the functions $\chi_n^{-1}(\frac{1}{\nu_n})$, \  
$3-\frac{4}{n+1}$ \  and \ $\frac{\sqrt{n-1}}{e}$ are shown in Fig.~\ref{fig:nev_ukl_theta_lower}.
These functions are defined for an integer argument,
but when constructing the graphs, we use continuous lines of various type.
The inequality 
 $$\chi_n^{-1}\left(\frac{1}{\nu_n}\right)> 3-\frac{4}{n+1}$$ 
 holds for $n\geq 53$.

\begin{table}[!htbp]
\begin{center}
\caption{\label{tab:nev_uhl_theta_lower2} Lower estimates of $\theta_n$ \ $(55\leq n \leq 109)$}
\bgroup
$
\def\arraystretch{1.44}
\begin{array}{|c|c|c|c||c|c|c|c|}
\hline
{\bf n} & {\bf \chi_n^{-1}\left(\frac{1}{\nu_n}\right) } & {\bf 3-\frac{4}{n+1}} & \text{\bf Max} &
{\bf n} & {\bf \chi_n^{-1}\left(\frac{1}{\nu_n}\right) } & {\bf 3-\frac{4}{n+1}} & \text{\bf Max} \\
\hline
 55 & 2.9669 & 2.9286 & 2.9669 & 82 & 3.5595 & 2.9518 & 3.5595 \\
 \hline
 56 & 3.0034 & 2.9298 & 3.0034 & 83 & 3.5546 & 2.9524 & 3.5546 \\
 \hline
 57 & 3.0315 & 2.931 & 3.0315 & 84 & 3.5899 & 2.9529 & 3.5899 \\
 \hline
 58 & 3.0583 & 2.9322 & 3.0583 & 85 & 3.6053 & 2.9535 & 3.6053 \\
 \hline
 59 & 3.058 & 2.9333 & 3.058 & 86 & 3.6355 & 2.954 & 3.6355 \\
 \hline
 60 & 3.0877 & 2.9344 & 3.0877 & 87 & 3.6307 & 2.9545 & 3.6307 \\
 \hline
 61 & 3.1173 & 2.9355 & 3.1173 & 88 & 3.667 & 2.9551 & 3.667 \\
 \hline
 62 & 3.1465 & 2.9365 & 3.1465 & 89 & 3.6803 & 2.9556 & 3.6803 \\
 \hline
 63 & 3.1463 & 2.9375 & 3.1463 & 90 & 3.7099 & 2.956 & 3.7099 \\
 \hline
 64 & 3.1849 & 2.9385 & 3.1849 & 91 & 3.7052 & 2.9565 & 3.7052 \\
 \hline
 65 & 3.2039 & 2.9394 & 3.2039 & 92 & 3.7401 & 2.957 & 3.7401 \\
 \hline
 66 & 3.2375 & 2.9403 & 3.2375 & 93 & 3.755 & 2.9574 & 3.755 \\
 \hline
 67 & 3.2322 & 2.9412 & 3.2322 & 94 & 3.783 & 2.9579 & 3.783 \\
 \hline
 68 & 3.2717 & 2.942 & 3.2717 & 95 & 3.7781 & 2.9583 & 3.7781 \\
 \hline
 69 & 3.2888 & 2.9429 & 3.2888 & 96 & 3.8121 & 2.9588 & 3.8121 \\
 \hline
 70 & 3.32 & 2.9437 & 3.32 & 97 & 3.8258 & 2.9592 & 3.8258 \\
 \hline
 71 & 3.3158 & 2.9444 & 3.3158 & 98 & 3.8542 & 2.9596 & 3.8542 \\
 \hline
 72 & 3.3529 & 2.9452 & 3.3529 & 99 & 3.8497 & 2.96 & 3.8497 \\
 \hline
 73 & 3.3703 & 2.9459 & 3.3703 & 100 & 3.8814 & 2.9604 & 3.8814 \\
 \hline
 74 & 3.4023 & 2.9467 & 3.4023 & 101 & 3.8964 & 2.9608 & 3.8964 \\
 \hline
 75 & 3.3973 & 2.9474 & 3.3973 & 102 & 3.9246 & 2.9612 & 3.9246 \\
 \hline
 76 & 3.4348 & 2.9481 & 3.4348 & 103 & 3.92 & 2.9615 & 3.92 \\
 \hline
 77 & 3.4521 & 2.9487 & 3.4521 & 105 & 3.9662 & 2.9623 & 3.9662 \\
 \hline
 78 & 3.4818 & 2.9494 & 3.4818 & 106 & 3.9936 & 2.9626 & 3.9936 \\
 \hline
 79 & 3.4769 & 2.95 & 3.4769 & 107 & 3.989 & 2.963 & 3.989 \\
 \hline
 80 & 3.5138 & 2.9506 & 3.5138 & 108 & 4.0841 & 2.9633 & 4.0841 \\
 \hline
 81 & 3.5288 & 2.9512 & 3.5288 & 109 & 4.034 & 2.9636 & 4.034 \\
 \hline
\end{array}
$
\egroup
\end{center}
\end{table}

\begin{table}[!htbp]
\begin{center}
\caption{\label{tab:nev_uhl_theta_lower3} Lower estimates of $\theta_n$ \ $(110\leq n \leq 118)$}
\bgroup
$
\def\arraystretch{1.44}
\begin{array}{|c|c|c|c||c|c|c|c|}
\hline
{\bf n} & {\bf \chi_n^{-1}(\frac{1}{\nu_n}) } & {\bf 3-\frac{4}{n+1}} & \text{\bf Max} &
{\bf n} & {\bf \chi_n^{-1}(\frac{1}{\nu_n}) } & {\bf 3-\frac{4}{n+1}} & \text{\bf Max} \\
\hline
 110 & 4.0615 & 2.964 & 4.0615 & 114 & 4.1281 & 2.9652 & 4.1281 \\
 \hline
 111 & 4.0568 & 2.9643 & 4.0568 & 115 & 4.1234 & 2.9655 & 4.1234 \\
 \hline
 112 & 4.0789 & 2.9646 & 4.0789 & 117 & 4.167 & 2.9661 & 4.167 \\
 \hline
 113 & 4.101 & 2.9649 & 4.101 & 118 & 4.1937 & 2.9664 & 4.1937 \\
 \hline
\end{array}
$
\egroup
\end{center}
\end{table} 

\section{Estimates of   Numbers $\xi_n^\prime$ and $\theta_n^\prime$}\label{nev_ukh_sec_5}

Let us give estimates of  numbers $\xi_n^\prime$ and $\theta_n^\prime$
  related to the situation
when  vertices of simplices located at vertices of the cube. By our definition,
 \begin{equation}\label{xi_teta_prime}
\xi_n^\prime:=\min\limits_{\ver(S)\subset \ver(Q_n)}\xi(S),  \qquad
 \theta_n^\prime:=\min\limits_{\ver(S)\subset \ver(Q_n)}\|P\|_{Q_n}.
 \end{equation}
In the right-hand equality, $S$ is the simplex having vertices at the nodes of a projector $P:C(Q_n)\to \Pi_1({\mathbb R}^n)$. 
  The sequence $ \theta^\prime_n$ was introduced by M.\,Nevskii in \cite{nevskii_mais_2003_10_1}.   He  proved
that  $\theta^\prime_n\geq \frac{1}{2e}\sqrt{n}$ , and if $ n + 1 $ is an Hadamard number, then 
$\theta^\prime_n\leq\sqrt{n+1}$.
More systematically
the numbers  \eqref{xi_teta_prime} were studied  by the authors  in  \cite{nev_ukh_saratov_2018}. 
Here we describe some results of this work.

First note general estimates for $\xi_n^\prime$ и $\theta_n^\prime$. Many of them are also estimates 
for $\xi_n$ and $\theta_n$.
The inequality holds
$$\theta^\prime_n>
\frac{1}{e}\sqrt{n-1}.$$
If $n\not=2$, then
$$\theta^\prime_n \leq \min \left( \frac{n+1}{2},
\frac{4\sqrt{e}}{3} \sqrt{n+2}+1 \right). 
$$
For each $n$, we have $\xi^\prime_n\geq n.$
If $n$ is even, then
$\xi^\prime_n>n.$ 
If $n+1$ is an Hadamard number, then $\xi^\prime_n=n$.
For $n>2$,
$$\xi^\prime_n \leq \frac{n^2-3}{n-1}.$$ 
Thus, $\theta^\prime_n\asymp n^{1/2}$,
$\xi^\prime_n\asymp n$. 


Assume $S$ is an $n$-dimensional     $(0,1)$-simplex having the maximum volume in $Q_n$ and   
$P$ is the corresponding interpolation projector. Then with constants not depending on  $n$ we have
$\|P\|_{Q_n}\asymp\theta^\prime_n$,
$\xi(S)\asymp \xi^\prime_n$. 

For $n\geq 3$, the sequence 
$\{\xi^\prime_n\}$ strictly increases. 


Consider the $n$-dimensional simplex $S^*$ for which the first $n$ vertices are the vertices of the cube
adjacent to $(1,\ldots,1)$ and the $(n+1)$th vertex coincides with $0$:
$$x^{(1)}=(0,1.\ldots,1,1), \  \ldots, \
x^{(n)}=(1,1.\ldots,1,0), \  x^{(n+1)}=(0,\ldots,0).$$  
The lengths of $n$ edges of  $S^*$ containing $x^{(n+1)}$ are equal to $\sqrt{n-1}$;
the lengths of the rest edges are $\sqrt{2}$. For $n=3$ (and only in this case)
$S^*$ is a regular simplex. 
Starting~$n=3$, the simplex $S^*$ has the following property
 (see \cite[Lemma 3.3]{hudelson_1996}): while replacing any vertex of
 $S^*$ by any point of $Q_n$, the volume of the simplex strictly decreases.
That is why for $n\geq 3$ the simplex   $S^*$ is called in 
\cite{hudelson_1996}  
{\it rigid.} This property holds also for  $n=1$. In  the two-dimensional case the volume of the simplex in such a way
does not increase.
If~and only if $1\leq n\leq 4$,
the volume of $S^*$ is the maximum possible for a simplex contained in $Q_n$.

 It was proved in  \cite{nevskii_monograph} that
$$\xi(S^*)=\left\{ \begin{array}{cc}
1,&n=1;\\  \\
4,&n=2;\\  \\
\frac{n^2-3}{n-1},&n\geq3. \\
\end{array}\right.
$$
 For each $n$, we have 
$d_1(S^*)=\ldots=d_n(S^*)=1,$ $\alpha(S^*)=n$. 
The one-dimensional case is the single when $\alpha(S^*)=
\xi(S^*)$.
This means that the simplex 
$\xi(S^*)S^*$ is circumscribed around $Q_n$ only for $n=1$.

 It turns out that for $1\leq n\leq 5$ the simplex $S^*$ delivers the exact values of
$\xi^\prime_n$, i.\,e., in these cases $\xi(S^*)=\xi^\prime_n$. 
The minimal $n$ satisfying
$\xi(S^*)>\xi^\prime_n$ is equal to $6$. 
Let $P^*:C(Q_n)\to \Pi_1({\mathbb R}^n)$ be the interpolation projector corresponding
to $S^*$.  The projector $P^*$ is extremal with respect to
 $\theta^\prime_n$, if $1\leq n\leq 4$. In  these dimensions, 
$\|P^*\|_{Q_n}=\theta^\prime_n$. The minimal $n$ with a property
$\|P^*\|>\theta^\prime_n$ equals $5$. Let us describe these cases 
in more details.

Suppose $n=1$.  We have $Q_1=S^*=[0,1]$, $\|P^*\|_{Q_1}=\xi(S^*)=1$,
hence,  $\theta_1=\xi_1=\theta^\prime_1=\xi^\prime_1=1$.
 
In the case $n=2$ there is only one (up to similarity) $(0,1)$-simplex, namely,
the simplex  $S^*$ with vertices $(0,1)$, $(1,0)$,  $(0,0)$.
It follows that $\theta^\prime_2=3$, $\xi^\prime_2=4$.

The case $n=3$ was analyzed by the first author (see, e.\,g., \cite{nevskii_monograph}).
The minima of  $\xi(S)$ and $\|P\|_{Q_3}$ are delivered by the simplex $S^*$ with vertices
$(0,1,1)$, $(1,0,1)$, $(1,1,0)$, $(0,0,0)$. 
We have $\xi_3=\xi^\prime_3=3$, $\theta_3=\theta^\prime_3=2$.
The same values of $\xi(S)$ and $\|P\|_{Q_3}$ are given also by the simplex
$S^{**}$ with vertices
$\left(\frac{1}{2}, 0, 0\right)$,
$\left(\frac{1}{2}, 1, 0\right)$,
$\left(0, \frac{1}{2}, 1\right)$,
$\left(1, \frac{1}{2}, 1\right)$
(and by similar simplices).
But  vertices of $S^{**}$ are not $(0,1)$-points.
There are no another extremal simplices in the three-dimensional case.

For $n=4$ and $n=5$, we have performed the full search of   $n$-dimensional simplices
satisfying the condition ${\rm ver}(S)\subset{\rm ver}(Q_n)$. It was discovered that 
$S^*$ and  $P^*$ give  both the minima  of 
$\xi(S)$ and $\|P\|_{Q_n}$. 
The following equalities hold: $\theta^\prime_4=\frac{7}{3}$,
 $\xi^\prime_4=\frac{13}{3}$.
The simplex $S^*$ remains an extremal also with respect to  $\xi_5$,
i.\,e.,~$\xi_5=\xi(S^*)=\frac{11}{2}$. It was obtained that 
$\theta_5$ $=$ $\frac{13}{5}$.  
The inequality $\theta^\prime_5\leq \frac{13}{5}$ was established  earlier \linebreak  in 
\cite[п.\,3.8.3]{nevskii_monograph}.

 Let us turn to the case $n=6$.
From consideration  $S^*$, it follows that 
$6\leq \xi^\prime_6 \leq \frac{33}{5}$.
The upper estimate may be improved.
We have made the full search of six-dimensional simplices satisfying the condition ${\rm ver}(S)\subset{\rm ver}(Q_6)$. 
For each simplex, one vertex was fixed and  the rest six vertices were chosen from the rest $2^6-1$ 
vertices of the cube. Thus, we were needed to look through 
$\binom{2^6-1}{6} = 67\,945\,521$ simplices.
Calculations were performed with the use of a special programme in  Wolfram Mathematica system (see, e.\,g., \cite{wolfram_2016}, \cite{klimov_uhalov_wolfram_2014}). 
The minimum of $\xi(S)$ equal to $\frac{25}{4}$ is delivered by the simplex with vertices
$(0, 0, 0, 0, 1, 1)$,
$(0, 0, 0, 1, 0, 1)$,
$(0, 1, 1, 1, 1, 0)$,
$(1, 0, 1, 1, 1, 0)$,
$(1, 1, 0, 1, 1, 0)$,
$(1, 1, 1, 0, 0, 1)$,
$(0, 0, 0, 0, 0, 0)$.
Therefore, $\xi^\prime_6=\frac{25}{4}$.
Simultaneously, we detected that the minimal $n$, for which
$S^*$ does not give an exact value of  $\xi^\prime_n$, is equal to~$6$.

From the results of  \cite[ п.\,3.8.3]{nevskii_monograph}, it follows that $\theta^\prime_6\leq 3$.
The full search of $(0,1)$-simplices in ${\mathbb R}^6$ have shown that the above simplex
minimizes also the norm of an interpolation projector. For the corresponding projector,
$\|P\|_{Q_6}=3$.
Thus, the exact value of  $\theta^\prime_6$ equals $3$.
Note that~$\|P^*\|_{Q_6}=\frac{17}{5}>\theta^\prime_6$. 

\begin{table}[h]
\begin{center}
\caption{ \label{nevsk_tab_ksitheta}
The numbers 
$\xi_n$, $\xi^\prime_n$, $\theta_n$, and  $\theta^\prime_n$ for $1\leq n\leq 7$} 

\medskip
\bgroup
\def\arraystretch{2.1}%
\begin{tabular}{|c|c|c|c|c|}
\hline
$n$ & $\xi_n$ & $\xi^\prime_n$ & $\theta_n$ &  $\theta^\prime_n$\\
\hline
 $1$ & $1$ &  $1$ & $1$ & $1$\\
 \hline
 $2$  & $\frac{3\sqrt{5}}{5}+1=2.34\ldots$  &   $4$ &  $\frac{2\sqrt{5}}{5}+1=1.89\ldots$   & $3$  \\
\hline
 $3$  & $3$  &   $3$  &    $2$   & $2$  \\
\hline
 $4$  & $4\leq \xi_4\leq \frac{19+5 \sqrt{13}}{9}=4.11\ldots$  & $\frac{13}{3}$ & 
 $\frac{11}{5} \leq \theta_4\leq \frac{3  (4+\sqrt{2})}{7}  =2.3203\ldots$    & $\frac{7}{3}$  \\
\hline
 $5$  & $5$   &   $\frac{11}{2}$   & $\frac{7}{3}\leq \theta_5< 2.448804$   & $\frac{13}{5}$ \\
\hline
 $6$  & $6\leq \xi_6< 6.0166$   &   $\frac{25}{4}$ &  $\frac{17}{7}\leq \theta_6\leq 2.60014$   &  3 \\
\hline
 $7$  & $7$   &  $7$   &    $\frac{5}{2}$   &  $\frac{5}{2}$  \\
 \hline
\end{tabular}
\egroup
\end{center}
\end{table}

Finally, consider the case $n=7$. Since  $8$ is an Hadamard number, there exists a seven-dimensional
regular simplex having  vertices at vertices of $Q_7$. This is, for example, the simplex
$ S $ with vertices
$(1,1,1,1,1,1,1)$, 
$(0,1,0,1,0,1,0)$,  
$(0,0,1,1,0,0,1)$,  
$(1,0,0,1,1,0,0)$,  
 $(0,0,0,0,1,1,1)$,  
$(1,0,1,0,0,1,0)$,  
$(1,1,0,0,0,0,1)$, 
 $(0,1,1,0,1,0,0)$.
We have $\xi(S)=7$.
From  inequalities $\xi^\prime_n\geq \xi_n\geq n$, it follows that $\xi^\prime_7 = \xi_7 =7$.
Now let us apply the estimate
\begin{equation}\label{nevsk_eq6}
\xi_n \leq
\frac{n+1}{2}\left( \theta_n-1\right)+1
\end{equation} 
which appears from the right-hand  inequality in \eqref{nev_ksi_P_ineq} for  $\Omega=Q_n.$
Since $\xi_7 =7$,  \eqref{nevsk_eq6} yields
$\theta_7 \geq \frac{5}{2}$. But  the projector corresponding to $S$ has the norm $\|P\|_{Q_7}=\frac{5}{2}$.
Therefore, $\theta^\prime_7 =\theta_7 = \frac{5}{2}$.

A collection of the known results related to  the numbers $\xi_n$, $\xi^\prime_n$, $\theta_n$, and  $\theta^\prime_n$
for~$1\leq n \leq 7$ is~given in Table  \ref{nevsk_tab_ksitheta}.

\section{Simplices Satisfying the Inclusions 
$S\subset Q_n\subset nS$}\label{nev_ukh_sec_6}


As is noted, the absorption index of the unit cube 
by a simplex contained  in~this cube is not less than the dimension of the space. This is equivalent to the inequality
$\xi_n\geq n$ which follows from
 \eqref{xi_S_alpha_S_geq_n}.
If for some $S\subset Q_n$ we have $\xi(S)=n,$  
then $\alpha(S)=\xi(S)=n$ and $\xi_n=n$. 
 Moreover, for each $i$,  holds $d_i(S)=1$.
On~the~contrary, the equality $\xi_n=n$ means that there exists a simplex
$S\subset Q_n$
with the property $\alpha(S)=\xi(S)=n$, i.\,e., $S\subset nS$.
However, simplices satisfying 
$S\subset Q_n\subset nS$
do exist not for any $n$.
The first such value is $n=2$.
Namely, $\xi_2=1+\frac{3\sqrt{5}}{5}=2.34 \ldots>2$; nontrivial proof is given in
 \cite{nevskii_monograph}. 
Below we note another way to see that $\xi_2>2$ .

Nowadays the equality  $\xi_n=n$ is proved for $n=5$, $n=9$ and for an infinite set of
$n$ such that $n+1$ is an Hadamard number. In particular, 
such are all odd $ n $ from an interval  $1\leq n\leq 11$.   No one
even $ n $ with the condition $\xi_n=n$  has not yet been found. Moreover, it remains unclear
does there exist at least one such even $ n $.

Simplices satisfying the inclusions
$S\subset Q_n\subset nS$, if exist, have a number of~remarkable properties. Such simplices were studied
by the authors in 
\cite{nev_ukh_mais_2017_24_5}.

Further $x^{(j)}$ are the vertices and
$\lambda_j(x)=
l_{1j}x_1+\ldots+
l_{nj}x_n+l_{n+1,j}$
are the basic Lagrange polynomials of $n$-dimensional simplex $S$. By $c(\Omega)$ we denote
the center of~gravity of a convex body  $\Omega$. By definition,
$$D_j:=\bigl \{x\in {\mathbb R}^n: \ -\frac{n-1}{n+1}\leq\lambda_j(x)\leq 1\bigr\}.$$
 
The following theorem combines the results of Theorem 1 and a number of corollaries from
\cite{nev_ukh_mais_2017_24_5}.

\smallskip 
{\bf Theorem 6.1.} {\it Suppose
$S\subset Q_n \subset nS$. Then the following relation hold:
\begin{equation}\label{max_lambda_j_eq_1}
\max\limits_{x\in \ver(Q_n)} \lambda_1(x)=
\ldots=
\max\limits_{x\in \ver(Q_n)} \lambda_{n+1}(x)=1, 
\end{equation}
\begin{equation}\label{max_minus_lambda_j_eq_nm1_np1}
\max\limits_{x\in \ver(Q_n)} \left(-\lambda_1(x)\right)=
\ldots=
\max\limits_{x\in \ver(Q_n)} \left(-\lambda_{n+1}(x)\right)=\frac{n-1}{n+1}, 
\end{equation}
\begin{equation}\label{theor_center}
c(S)=c(Q_n)=
\left(\frac{1}{2},\ldots,\frac{1}{2}\right),
\end{equation}
$$\ldots\subset \frac{1}{n^2}\,S\subset\frac{1}{n}\,Q_n
\subset S\subset Q_n\subset nS \subset n^2Q_n\subset n^3 S\subset\ldots,$$
$$\sum_{j=1}^{n+1} |l_{ij}|=2 \qquad(1\leq i\leq n), 
\sum_{i=1}^{n} |l_{ij}|=\frac{2n}{n+1} \qquad(1\leq j\leq n+1),$$
$$\sum_{i\leq n: \ l_{ij}\geq 0} l_{ij}=1-l_{n+1,j}, \
\sum_{i\leq n: \ l_{ij}< 0} |l_{ij}|=\frac{n-1}{n+1}+l_{n+1,j} \quad(1\leq j\leq n+1),$$
$$Q_n\subset\bigcap\limits_{j=1}^{n+1} D_j.$$ 

\smallskip The equalities similar to $(\ref{max_lambda_j_eq_1})$--$(\ref{max_minus_lambda_j_eq_nm1_np1})$
hold also in the case when maxima are taken over $x\in Q_n$.
}

\smallskip
Applying  \eqref{theor_center}, it is rather easy to prove that 
$\xi_2>2.$ Earlier the exact value of $\xi_2$ was obtained by M. Nevskii 
 (see  \cite{nevskii_monograph}). Let us note some other corollaries 
 of Theorem 6.1.

\smallskip
{\it Let $n$ be an even number. Then there does not exist  a simplex satisfying simultaneously
the condition
$S\subset Q \subset nS$ and $\ver(S)\subset \ver(Q_n)$.}

\smallskip
Assume the converse: there exists a simplex $S$ such that the vertices $x^{(j)}$
coincide with vertices of  $Q_n$ and
$S\subset Q_n\subset nS$.  By \eqref{theor_center}, 
$$c(S)=\frac{1}{n+1}\sum_{j=1}^{n+1} x^{(j)}=
\left(\frac{1}{2},\ldots,\frac{1}{2}\right),$$
hence,
$$u:=\sum_{j=1}^{n+1} x^{(j)}=\left(\frac{n+1}{2},\ldots,\frac{n+1}{2}\right).$$
Since $x^{(j)}\in \ver(Q_n)$, 
coordinates of $u$ are integer.
If  $n$ is even, the latter equality cannot be valid. Thus, we have a contradiction.

\smallskip {\it Suppose
$S\subset Q_n\subset nS$. Then this simplex $S$ has the following property: 
while replacing any vertex of $S$ by an arbitrary point of $Q_n$,  the volume of the simplex does not increase.}

\smallskip
Let us prove it.
Assume $1\leq j\leq n+1$ and $y\in Q_n$. Since $$-\frac{n-1}{n+1}\leq\lambda_j(y)\leq 1, \quad
\lambda_j\left(x^{(j)}\right)=1,$$ 
we have
\begin{equation}\label{dist_leq_dist}
{\mathrm {dist}}(y; B_j)\leq {\mathrm {dist}}\left (x^{(j)}; B_j\right). 
\end{equation}
Here $B_j$ is the $(n-1)$-dimensional hyperplane given by the equation $\lambda_j(x)=0$ and
$${\mathrm {dist}} (y; B_j):=\min_{z\in B_j} \|y-z\|$$
is the distance from $y$ to $B_j$. Denote by $S^\prime$ a simplex 
with vertices
$x^{(1)},$ $\ldots,$ $x^{(j-1)},$ $y,$ $x^{(j+1)},$ $\ldots,$ $x^{(n+1)}$, i.\,e.,
the convex hull of the points $x^{(k)}$ for $k\ne j$ and the point~$y$. From (\ref{dist_leq_dist}), it follows that
$\vo \left(S^\prime\right)\leq \vo (S)$.

By the terminology of  \cite{hudelson_1996}, a simplex $S\subset Q_n$,
for which replacing any vertex by any point of $ Q_n $ decreases its volume, is called
{\it rigid.} According to the   above result, a simplex $S$ satisfying the inclusions
$S\subset Q_n\subset nS$, may be called {\it quasi-rigid}.


\section{Perfect Simplices}\label{nev_ukh_sec_7}


This section is  based on   the authors'  paper \cite{nev_ukh_beitrage_2018}, see also \cite{nev_ukh_posobie_2020}.  

By our definition, {\it $S$ is a perfect simplex with respect to a cube $Q$} if  $S\subset Q\subset \xi_n S$ and the cube $Q$ is inscribed into the simplex $\xi_n S$, i.\,e., the boundary of
$\xi_n S$ contains all the vertices of
$Q$. A simplex perfect with respect to the unit cube $Q_n$ is called simply {\it perfect}.
Today we know the only three dimensions $n$ when there exist perfect simplices in  ${\mathbb R}^n$.
These are $n=1$, $n=2$, and $n=5$; in all the cases   $\xi_n=n$.
 
Consider the question on existence of perfect simplices in  ${\mathbb R}^n$ for various $n$.
We need the following proposition proved in \cite{nev_ukh_beitrage_2018}.

\smallskip
{\bf Theorem 7.1.} {\it
Suppose $S\subset Q_n$. A vertex $v$ of $Q_n$ belongs to the $j$th face of the simplex  
$\xi(S)S$  when and only when 
\begin{equation}\label{nev_uhl_cub_lambd}
-\lambda_j(v)=\max_{1\leq k\leq n+1, x\in \ver(Q_n)} (-\lambda_k(x)). 
\end{equation}
}

\smallskip
By the $j$th face of the simplex $\xi(S)S$, we mean the $(n-1)$-dimensional face  contained in the hyperplane
$\lambda_j(x)={\rm const}$. As above, $\lambda_j$ is the $j$th basic Lagrange
polynomial related to $S$.

\medskip
1) $n=1$.

This case is very simple. For the segment $S=[0,1]$,  we have $S= Q_1$.
Hence, $\xi_1= 1$  and $S$ is a unique perfect simplex.
The equality $\xi_1= 1$ also follows from Theorem ~2.1 since  
2 is an Hadamard number.

\smallskip
2) $n=2$.

There are no perfect simplices in ${\mathbb R}^2$.  A unique (up to rotations) simplex $S$
satisfying $\xi(S)=\xi_2$ is described  in
\cite[\S\,2.4]{nevskii_monograph}.  Obviously, $S$ is not perfect, since the zero-vertex of  $Q_2$
is not contained in the boundary of $\xi(S)S$.
\smallskip

3) $n=3$.

The exact value $\xi_3=3$ is given in \cite[\S\,2.4]{nevskii_monograph}.
The equality $\xi(S)=3$  holds for two (up to equivalence) tetrahedrons.  Let us consider these cases.

Since for $n=3$  the conditions of Theorem 2.1 are fulfilled, we have $\xi_3=3$. There exists a regular
simplex  $S$ inscribed into  $Q_3$ and satisfying $Q_3\subset 3S$.
We may take the simplex $S_1$ with vertices
$(0, 0, 0)$,
$(1, 1, 0)$,
$(1, 0, 1)$,
$(0, 1, 1)$.
For this simplex, 
$$
{\bf A}=
\left(
\begin{array}{cccc}
 0 & 0 & 0 & 1 \\
 1 & 1 & 0 & 1 \\
 1 & 0 & 1 & 1 \\
 0 & 1 & 1 & 1 \\
\end{array}
\right),\quad
{\bf A}^{-1}=\frac{1}{2}
\left(
\begin{array}{rrrr}
 -1 & 1 & 1 & -1 \\
 -1 & 1 & -1 & 1 \\
 -1 & -1 & 1 & 1 \\
 2 & 0 & 0 & 0 \\
\end{array}
\right).
$$
Formulae \eqref{alpha_d_i_formula}--\eqref{d_i_l_ij_formula} give 
$d_1(S_1)=d_2(S_1)=d_3(S_1)=1$, $\alpha(S_1)= 3$.
The basic Lagrange polynomials are
\begin{align*}
\lambda_1(x)&=  \frac{1}{2} \left(-x_1-x_2-x_3+2\right), &
\lambda_2(x)&= \frac{1}{2} \left(x_1+x_2-x_3\right), \\
\lambda_3(x) &= \frac{1}{2} \left(x_1-x_2+x_3\right),  &
 \lambda_4(x) &= \frac{1}{2} \left(-x_1+x_2+x_3\right).
\end{align*}
Calculation shows that 
\begin{equation}\label{nev_uhl_max12}
\xi(S_1)=4\max_{1\leq k\leq 4} \,\max_{x\in \ver(Q_3)}(-\lambda_k(x))+1=
4 \cdot  \frac{1}{2} +1 =3.
\end{equation}
 For each  $k$, maximum in \eqref{nev_uhl_max12}  is delivered by the only one vertex of $Q_3$:
\begin{equation}\label{nev_uhl_lam3vertcond}
-\lambda_1(1,1,1) = -\lambda_2(0,0,1) = -\lambda_3(0,1,0)= -\lambda_4(1,0,0) = \frac{1}{2}.
\end{equation}
By Theorem~7.1 and (\ref{nev_uhl_lam3vertcond}), the extremal vertices of the cube,
i.\,e., the vertices
$(1,1,1)$, $(0,0,1)$, $(0,1,0)$, $(1,0,0)$, belong to the faces of the simplex 
 $3S_1$. Moreover, each face contains a unique vertex. The rest four vertices of the cube are inside the
 simplex  $3S_1$.
Therefore, $S_1$, though satisfies the inclusions 
$S \subset Q_3 \subset 3S_1$, is not a perfect simplex.
Largest segments in  $S_1$ parallel to the axes intersect at the center of $Q_3$.

Now consider the simplex $S_2$ with vertices
$\left(\frac{1}{2},0,0\right)$,  
$\left(\frac{1}{2},1,0\right)$, 
$\left(0,\frac{1}{2},1\right)$,\linebreak
$\left(1,\frac{1}{2},1\right)$.
For $S_2$, 
$$
{\bf A}=
\left(
\arraycolsep=5pt\def\arraystretch{1.2}
\begin{array}{cccc}
 \frac{1}{2} & 0 & 0 & 1 \\
 \frac{1}{2} & 1 & 0 & 1 \\
 0 & \frac{1}{2} & 1 & 1 \\
 1 & \frac{1}{2} & 1 & 1 \\
\end{array}
\right),\quad
{\bf A}^{-1}=\frac{1}{2}
\left(
\arraycolsep=5pt\def\arraystretch{1.2}
\begin{array}{rrrr}
 0 & 0 & -1 & 1 \\
 -1 & 1 & 0 & 0 \\
 -\frac{1}{2} & -\frac{1}{2} & \frac{1}{2} & \frac{1}{2} \\
 1 & 0 & \frac{1}{2} & -\frac{1}{2} \\
\end{array}
\right).
$$
By our formulae,
$d_1(S_2)=d_2(S_2)=d_3(S_2)=1$, $\alpha(S_2)= 3$. We have
\begin{align*}
\lambda_1(x)&=  -x_2-\frac{1}{2}\, x_3+1, &
\lambda_2(x)&= x_2-\frac{1}{2}\, x_3, \\
\lambda_3(x) &= \frac{1}{2} \left(-2 x_1+x_3+1\right),  &
 \lambda_4(x) &= \frac{1}{2} \left(2 x_1+x_3-1\right),
\end{align*}
\begin{equation} \label{nev_uhl_max12vert}
\begin{array}{l}
- \lambda_1(0,1,1) = - \lambda_1(1,1,1) = - \lambda_2(0,0,1) = - \lambda_2(1,0,1) = \\ \\
- \lambda_3(1,0,0) = - \lambda_3(1,1,0) = - \lambda_4(0,0,0) = - \lambda_4(0,1,0) = \frac{1}{2}.
\end{array}
\end{equation}
Hence,
$$
\xi(S_2)=4\max_{1\leq k\leq 4} \,\max_{x\in \ver(Q_3)}(-\lambda_k(x))+1=4 \cdot  \frac{1}{2} +1 =3.
$$
In this case, any vertex of the cube is an extremal.
By Theorem~7.1 and  (\ref{nev_uhl_max12vert}), all the vertices of  $Q_3$ belong
to the faces of $3S_2$. 
This means that $S_2$, differently from $S_1$, is a perfect simplex.

\smallskip

4) $n=5$.

 As we  marked, $\xi_5=5$. Let us consider the following family of five-dimensional simplices $V=V(s,t)$. This
 family was introduced and studied by the authors in~\cite{nev_ukh_beitrage_2018}.
The vertices of
$V$ are the points
$\left(s,1,\frac{1}{3},1,1\right)$, 
$\left(s, 0,  \frac{1}{3}, 1, 1\right)$, $\left(s,2-3 t,\frac{1}{3},0,1\right)$,
$\left(2-3 s,t,0,\frac{1}{3},0\right)$, $\left(0, t, 1, \frac{1}{3}, 0\right)$, 
$\left(1, t, 1, \frac{1}{3}, 0\right)$. The vertices matrix is
$$
{\bf A}={\bf A}(s,t)=
\left(
\arraycolsep=3pt\def\arraystretch{1.3}
\begin{array}{cccccc}
 s & 1 & \frac{1}{3} & 1 & 1 & 1 \\
 s & 0 & \frac{1}{3} & 1 & 1 & 1 \\
 s & 2-3 t & \frac{1}{3} & 0 & 1 & 1 \\
 2-3 s & t & 0 & \frac{1}{3} & 0 & 1 \\
 0 & t & 1 & \frac{1}{3} & 0 & 1 \\
 1 & t & 1 & \frac{1}{3} & 0 & 1 \\
\end{array}
\right).
$$
For any $s,t\in {\mathbb R}$, the simplex  $V$ is nondegenerate and $\vo(V)=\frac{1}{120}.$ 
The center of gravity of the simplex coincides with the center of the cube.
The matrix ${\bf A}^{-1}$  has the form
$$
{\bf A}^{-1}={\bf A}^{-1}(s,t)=
\left(
\arraycolsep=5pt\def\arraystretch{1.4}
\begin{array}{cccccc}
 0 & 0 & 0 & 0 & -1 & 1 \\
 1 & -1 & 0 & 0 & 0 & 0 \\
 0 & 0 & 0 & -1 & 3 s-1 & 2-3 s \\
 2-3 t & 3 t-1 & -1 & 0 & 0 & 0 \\
 3 t-\frac{4}{3} & \frac{5}{3}-3 t & \frac{2}{3} & -\frac{2}{3} & 3 s-\frac{5}{3} & \frac{4}{3}-3 s \\
 -\frac{2}{3} & \frac{1}{3} & \frac{1}{3} & 1 & 2-3 s & 3 s-2 \\
\end{array}
\right).
$$
The inclusion  $V\subset Q_5$ takes places if and only if 
$\frac{1}{3}\leq s \leq \frac{2}{3},$ $\frac{1}{3}\leq t \leq \frac{2}{3}.$
Calculations described in \cite{nev_ukh_beitrage_2018} lead to the following conclusions.

Under the condition
\begin{equation}\label{nev_uhl_xi_is_5_cond}
 s\in\left[\frac{4}{9},\frac{5}{9}\right], \quad  t \in\left[\frac{4}{9},\frac{5}{9}\right],
\end{equation}
we have the equalities $d_i(V)=1,$ $\alpha(V)=\xi(V)=5.$
Moreover, applying Theorem 7.1, we see that in this case all the vertices of the cube
$Q_5$ belong to the boundary of the simplex  5V. This means that  {\it under the key condition  \eqref{nev_uhl_xi_is_5_cond} each simplex $V(s,t)$
is perfect.} 

Note that if
$s\notin\left[\frac{4}{9},\frac{5}{9}\right]$ or $t\notin\left[\frac{4}{9},\frac{5}{9}\right]$, then  $\xi(V(s,t))>5$. 

 The fact that there are perfect simplices in ${\mathbb R}^5$ is the basic result of 
 the authors' paper \cite{nev_ukh_beitrage_2018}.  It is still unknown whether there exist perfect simplices for odd   $n>5$. Also it remains unclear whether
 there exist perfect simplices for  at least one even $ n $.

\section{ Equisecting Simplices}\label{nev_ukh_sec_8}

The notion {\it  equisection property} was introduced by the authors in  \cite{nev_ukh_mais_2017_24_5} and also was studied by E.\,Esipova and A.\,Ukhalov in  \cite{esipova_2017}--\cite{esipova_ukhalov_2018}.

We call $S\subset Q_n$ {\it an equisecting simplex (or a simplex having the equisection property)}, if all the $(n-1)$-dimensional hyperplanes containing the
faces of  $S$ cut off from the cube (outward from the simplex)  the parts of equal volumes.  
A simplex $S\subset Q_n$ is called here {\it an extremal} if $\xi(S)=\xi_n$.

The equations of $(n-1)$-dimensional hyperplanes containing
 the
faces of  $S$ have the form  
$\lambda_j(x)=0$. The $j$th hyperplane contains all vertices of the simplex except the $j$th one.
 Denote by $v_j$ the volume of a convex  body cut off from the cube by the $j$th hyperplane and not containing the simplex. Thus,
$$
v_j = \operatorname{vol}(\,Q_n \cap \{x\in {\mathbb R}^n  :\ \lambda_j(x)\leq 0\}\,).
$$ 
The simplex is an equisecting if and only if $v_1 = v_2=\ldots =v_{n+1}$. 

We have noticed that some  extremal simplices are also equisecting.  
Obviously, the extremal one-dimensional simplex $S=[0,1]$ has this property, since the faces of $S$
cut off from  $Q_1=[0,1]$ the zero-length segments. The triangle with vertices
$(0,0),$ $(1,\tau)$, $(\tau,1)$, where $\tau=\frac{3-\sqrt{5}}{2}=0.38196601\ldots$,
is an equisecting simplex on plane.
M.\,Nevskii showed (see~\cite{nevskii_monograph}) that this is a unique (up to rotations)
extremal simplex for $n=2$.
The area of any triangle cut off from  $Q_2$ equals
$\frac{3-\sqrt{5}}{4}=0.1909\ldots$
It is easy to see that both extremal simplices for $n=3$ described in Section 7 also are equisecting.
Therefore, for $1\leq n \leq 3$, all the extremal simplices  have the equsection property.

Note that the equisection property does not imply an extremality. One can construct
equisecting simplices being not extremal. A problem of connection between equisection and extremality
is rather interesting.
It was discovered that the authors'  first conjecture that all the extremal simplices are also equisecting
is not true
(see~\cite{nev_ukh_mais_2017_24_5}). 
In particular, for 
$n=5$, $n=7$, and $n=9$, A.\,Ukhalov has built extremal simplices without the equisection property.
Yet it was discovered that for all extremal simplices the volumes $v_j$ obey some regularities. Namely, the set of
these  values  is always limited to only two quantities. 
Further note  that equisecting simplices do exist also for $n>3$. In \cite{esipova_ukhalov_2018},
it was proved that the five-dimensional extremal  simplex
$V=V(s,t)$ described in Section 7  has the equisection property. Let us give the corresponding
formulae from \cite{esipova_ukhalov_2018}  for the cut-off volumes: 

$$
v_1(t) = 
\left\{
\arraycolsep=5pt\def\arraystretch{1.4}
\begin{array}{ll}
 \frac{81 t-35}{1458 t^2-1620 t+432}+\frac{1}{2}, & -\infty <t<\frac{1}{3}, \\
 -\frac{t}{2}+\frac{2}{54-81 t}+\frac{4}{9}, &  \frac{1}{3}<t<\frac{4}{9}, \\
 \frac{1}{3}, & \frac{4}{9}<t<\frac{2}{3}, \\
\frac{t}{2}+\frac{2}{81 t-36}-\frac{1}{9}, & \frac{2}{3}<t<\frac{7}{9}, \\
\frac{55-81 t}{1458 t^2-1620 t+432}+\frac{1}{2}, & \frac{7}{9}<t<+\infty,
\end{array}
\right.
$$
$$
v_2(t) =
\left\{
\arraycolsep=5pt\def\arraystretch{1.4}
\begin{array}{ll}
\frac{81 t-26}{1458 t^2-1296 t+270}+\frac{1}{2}, & -\infty <t<\frac{2}{9}, \\
-\frac{t}{2}+\frac{2}{45-81 t}+\frac{7}{18}, & \frac{2}{9}<t<\frac{1}{3}, \\
 \frac{1}{3}, & \frac{1}{3}<t<\frac{5}{9}, \\
\frac{81 t^2-36 t+7}{162 t-54}, & \frac{5}{9}<t<\frac{2}{3}, \\
\frac{46-81 t}{1458 t^2-1296 t+270}+\frac{1}{2}, & \frac{2}{3}<t<+\infty,
\end{array}
\right.
$$
\vspace{3mm}

$$
v_3 =  \frac{1}{3}, \qquad v_4 =  \frac{1}{3}, \qquad
v_5(s) =  v_2(s), \qquad v_6(s) = v_1(s).
$$

The expressions for $v_5(s)$ and $v_6(s)$ are the same as for  $v_2(t)$ and $v_1(t)$ 
relatively. An analysis of the expressions for $v_1(t)$ and $v_2(t)$ shows that
  $v_2(t) = v_1\left(t+\frac{1}{9}\right)$.
Therefore, the behaviour of $v_j$ for $j=1,2,5, 6$ is described, in fact, by the same function.
 
\begin{figure}[h]
\centering
\includegraphics[width=1\textwidth]{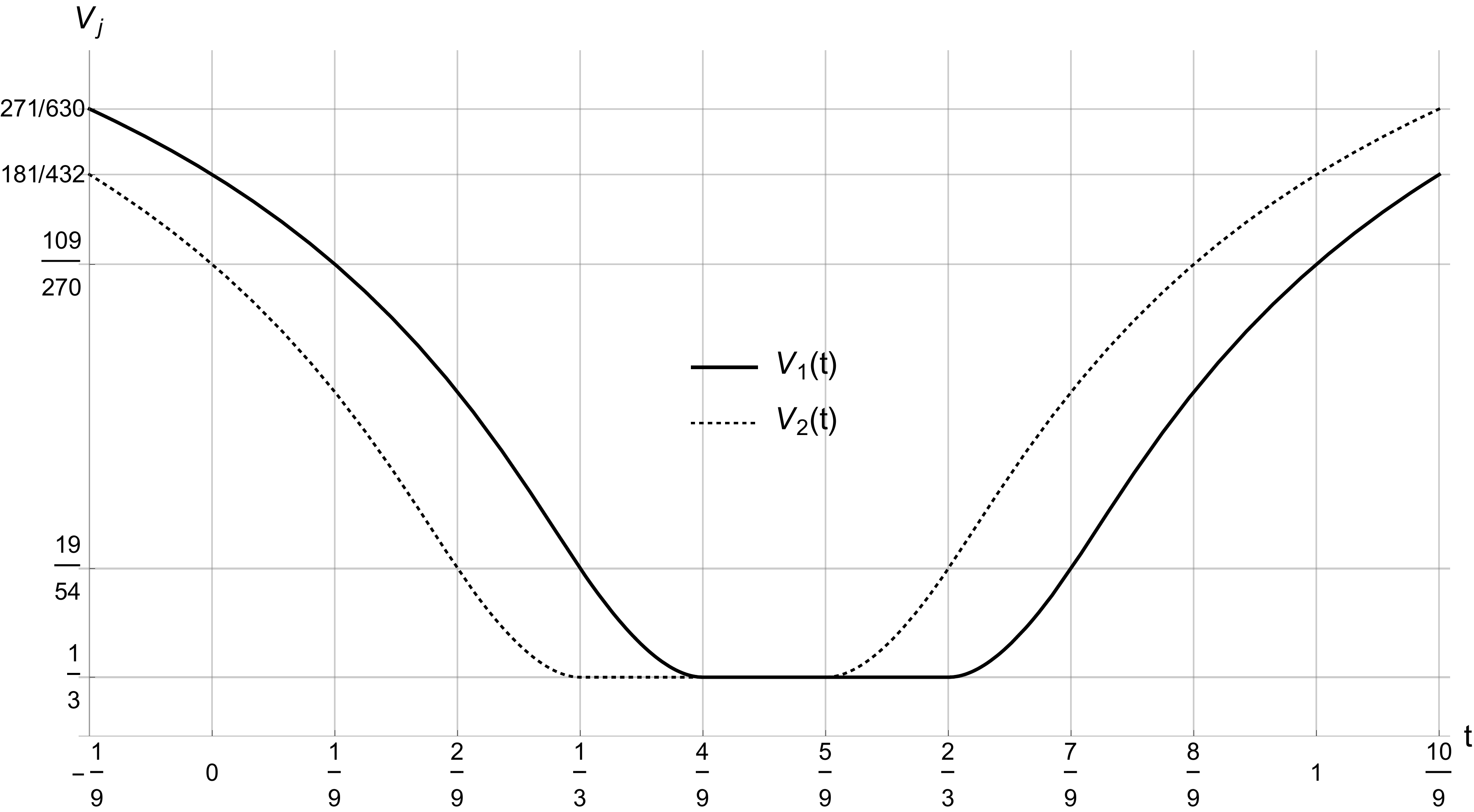}
\caption{Graphs of the functions $v_1(t)$ and $v_2(t)$}
\label{fig:esip_uhl_perf_v1v2}
\end{figure}

In Fig.~\ref{fig:esip_uhl_perf_v1v2}, we present  graphs of the functions $v_1(t)$ and $v_2(t)$ on
the segment $\left[-\frac{1}{9},\frac{10}{9}\right]$. 

By the above formulae,  if
$s,t\in\left[\frac{4}{9},\frac{5}{9}\right]$,
then $v_j=\frac{1}{3}$, $j=1,\ldots,6$.
It is for those $ s $ and $ t $ for which $ V (s, t) $ is perfect, this simplex is also equisecting.

We see that all  known today  perfect simplices are also equisecting.  
This allows the authors to conjecture the following: {\it all the perfect simplices have the equisection property.}

\section{ Properties of $(0,1)$-Matrices of Order $n$ \\ Having  Maximal Determinant }\label{nev_ukh_sec_9}

In this section, we describe the results presented in  \cite{nev_ukh_matzam SVFU_2019}. 
Complementary materials to this paper are given in \cite{nev_ukh_dataset_2019}.

As we mark in Section 1,  if a simplex
$S\subset Q_n$ has the maximum volume, then the axial diameters of $S$ are equal to $1$.
For this simplex $\alpha(S)=n$. By
\eqref{alpha_d_i_formula}, this is equivalent to the equalities
 $d_i(S)=1$. These are the basic facts for our analysis.

Denote by $h_n$ and $g_n$  maximal determinants of
$(0,1)$- and~$(-1,1)$-matrices of order $n$ respectively.
Define $\nu_n$  as the maximum volume of an  $n$-dimensional simplex
contained in
$Q_n$. We have
$g_{n+1}=2^{n}\,h_n$,
$h_n=n!\,\nu_n$ (see~\cite{hudelson_1996}). 

If a $(0,1)$-matrix of order $n$ having maximal determinant  is known, it is possible to 
construct a maximal volume simplex in $Q_n$.
Let us enlarge the row set of such a determinant by
the row $(0,\ldots,0).$ Then simplex $S$ with these vertices 
is contained in $Q_n$ and has the maximal possible volume.
Indeed, consider for $S$ the vertices matrix ${\bf A}$.
Then 
$$\vo(S)=\frac{|\det({\bf A})|}{n!}=\frac{h_n}{n!}=\nu_n.$$
One can also obtain nonzero vertices of an $n$-dimensional simplex  in $Q_n$
with maximal volume in the same way using 
the columns of a $(0,1)$-matrix with maximal determinant of order $n$. 
On the other hand, if an $n$-dimensional simplex $S\subset Q_n$ 
with the $0$-vertex
has the maximal possible volume, 
then the coordinates of its  nonzero vertices, 
being written in rows or in columns, form $(0,1)$-matrices with maximal determinant 
of order $n$.

{\bf Theorem 9.1.}
{\it  Suppose {\bf M} is an arbitrary nondegenerate $(0,1)$-matrix of order $n$.
Let $\bf A$ be the matrix of order $n+1$ which is obtained
from ${\bf M}$
by adding the $(n+1)$th row $(0,0,\ldots,0,1)$ and the $(n+1)$th column
consisting of $1$'s. Denote ${\bf A}^{-1}=(l_{i,j}).$ Then for each $i=1,\ldots,n$
\begin{equation}\label{nev_ukh_sum_abs_l_ij_geq2} 
\sum_{j=1}^{n+1} |l_{i,j}|\geq 2. 
\end{equation}
If $|\det({\bf M})|=h_n$, then 
for all $i=1,\ldots,n$ 
\begin{equation}\label{nev_ukh_sum_abs_l_ij_eq2}
 \sum_{j=1}^{n+1} |l_{i,j}|= 2. 
\end{equation}
If for some $i$ we have  the strong inequality
\begin{equation}\label{nev_ukh_sum_abs_l_ij_greater2}
 \sum_{j=1}^{n+1} |l_{i,j}|> 2,
\end{equation}
then  $|\det({\bf M})|<h_n$. 
}

\smallskip
{\it Proof.} Consider the simplex $S$ with the $0$-vertex and the rest
vertices coinciding in the coordinate
form with the rows of ${\bf M}$. Obviously, 
${\bf A}$
is the vertices matrix of $S$. Since $\det({\bf A})=\det({\bf M})\ne 0,$
the simplex $S$ is nondegenerate. 
The inclusion 
$S\subset Q_n$ means that $d_i(S)\leq 1$. 
Being combined with  \eqref{d_i_l_ij_formula},
this gives 
(\ref{nev_ukh_sum_abs_l_ij_geq2}).

Now assume  $|\det({\bf M})|=h_n$. Then we have
$$\vo(S)=\frac{|\det({\bf A})|}{n!}=\frac{|\det({\bf M})|}{n!}=\frac{h_n}{n!}=\nu_n.$$
Consequently, $S$ is a simplex in $Q_n$ 
with the maximal possible volume. 
 All the axial diameters $d_i(S)$ are equal to 1.
Applying  \eqref{d_i_l_ij_formula}, 
we obtain the conditions 
(\ref{nev_ukh_sum_abs_l_ij_eq2}).
Finally, if for some $i$ the strong inequality
(\ref{nev_ukh_sum_abs_l_ij_greater2}) holds then $d_i(S)<1.$  
In this case, the volume of $S$ is not maximal
and $|\det({\bf M})|<h_n$. 
\hfill$\Box$

\smallskip
Note that a $(0,1)$-matrix  $\bf T$  of order $n$  having maximal determinant $h_n$ can be obtained from 
a $(-1,1)$-matrix $\bf U$  of order $n+1$  with maximal determinant $g_{n+1}$.
It is sufficient to apply to ${\bf U}$ a special procedure described in  \cite{nev_ukh_mais_2018_25_1}. 
If at the exit will turn out $|\det({\bf T})|<h_n,$ then at the entrance we had $|\det({\bf U})|<g_{n+1}.$  
Using the above theorem, in some
cases, it is possible to show that the corresponding determinants are non-maximal.

\begin{figure}[!htbp]
\begin{center}
\includegraphics[width=\textwidth]{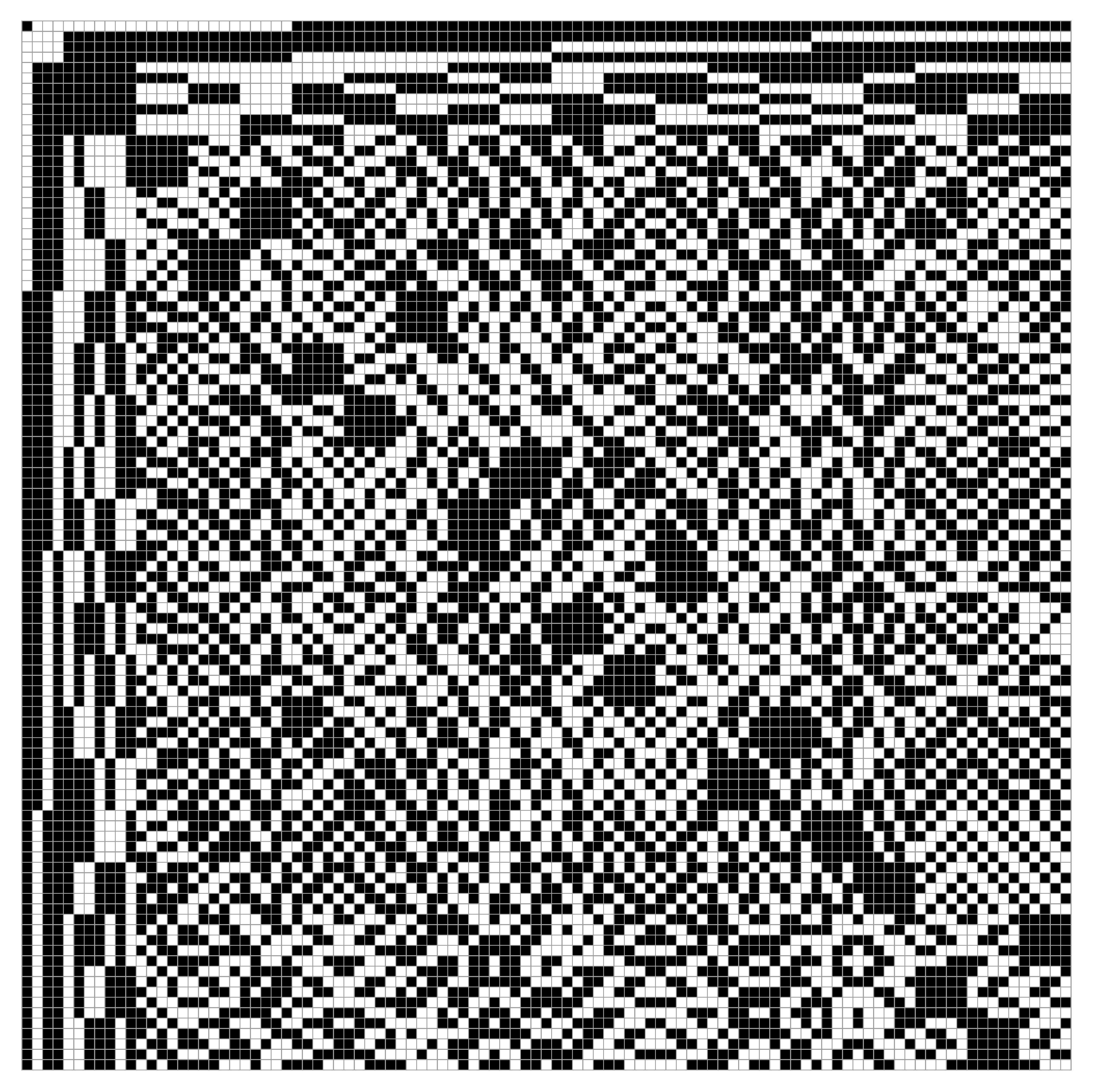}
\end{center}
\caption{$(-1,1)$-matrix matrix ${\bf U}$ of order $101$. White squares denote $-1$'s, black squares denote~$1$'s}
\label{fig:nev_uhl_matrix_plus_minus_ones}
\end{figure}

 \begin{figure}[!htbp]
\begin{center}
\includegraphics[width=\textwidth]{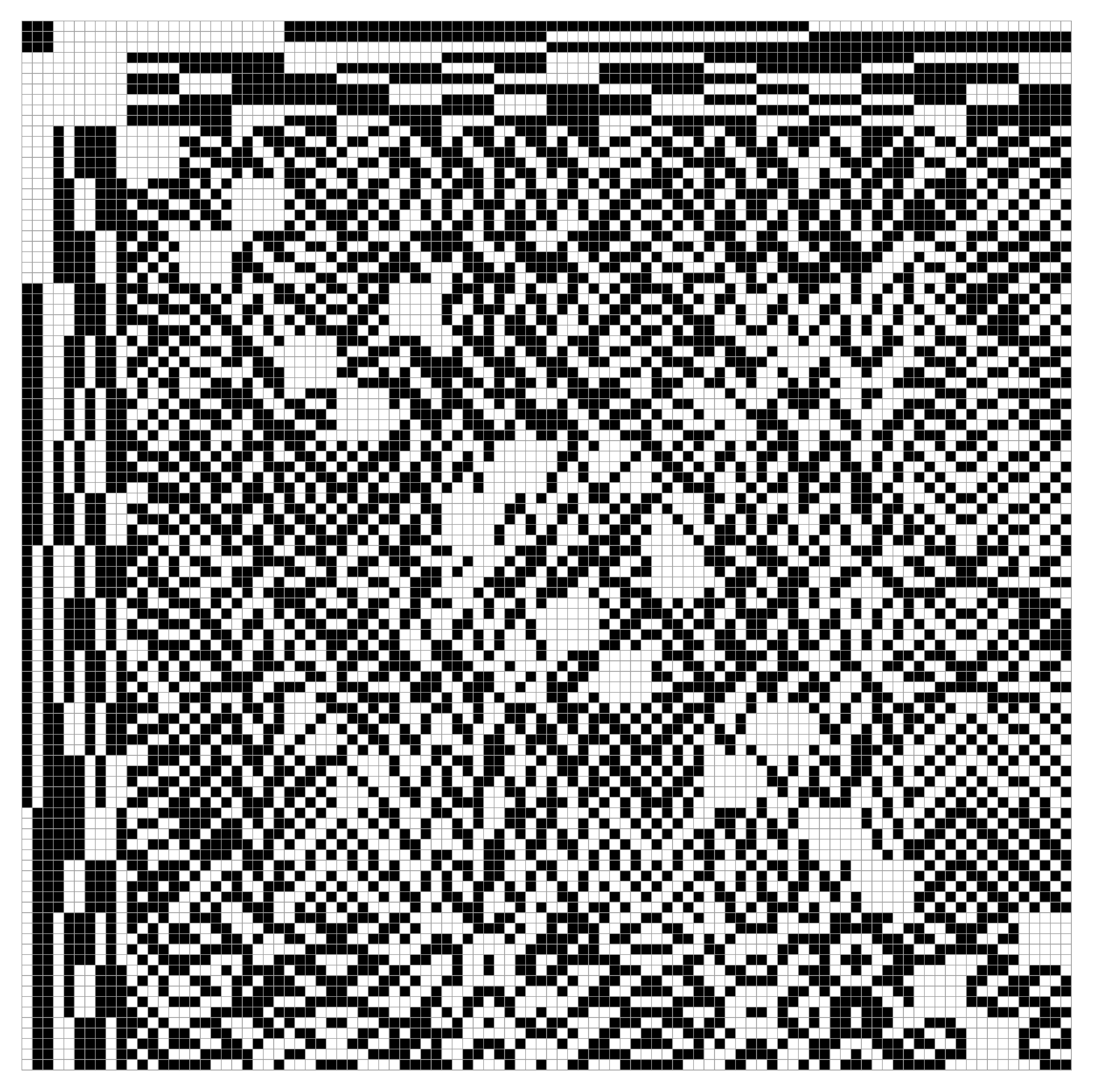}
\end{center}
\caption{$(0,1)$-matrix ${\bf T}$ of order $100$. White squares denote $0$'s, black squares denote $1$'s}
\label{fig:nev_uhl_matrix_zero_ones}
\end{figure}

A collection of matrices having maximal known determinants of orders $1,$ $\dots,$ $119$ 
can be found on the site \cite{Nevskii_bib_7}.
 Consider the matrix of order 101   constructed in~2003 by W.\,Orrick and B.\,Solomon.  
They don't claim that the given $(-1,1)$-matrix 
 has the largest possible determinant but 
 only state that its determinant  <<surpasses the previous record>>.
 In graphical form, the original $(-1,1)$-matrix $\bf U$ 
of order 101 and the resulting $(0,1)$-matrix
$\bf T$  of order 100 are shown 
in Fig. \ref{fig:nev_uhl_matrix_plus_minus_ones}--\ref{fig:nev_uhl_matrix_zero_ones}.

 Our  calculations performed by the use of the  Wolfram Mathematica system
show that a huge in absolute value 81-digit determinant of $\bf T$ :
$$|\det({\bf T})|=
2\cdot3^2\cdot5^{97}\cdot79=$$
$$89740816577942302983638241586569761487623964058002457022666931152343750,$$
\smallskip
yet is not maximal.
Applying Theorem 9.1 to   ${\bf M}={\bf T}$, we have 
$$
\sum_{j=1}^{101} |l_{ij}| = 
\left\{
 \arraycolsep=3pt\def\arraystretch{1.2} 
\begin{array}{cl}
 \frac{1438}{711}=2.0225\ldots, & i=1,2,3, \\
 \frac{490}{237}=2.0675\ldots, & i=4, \\
 2, & i=5,\ldots,100. \\
\end{array}
\right. 
$$
For $i\leq 4$, the strict inequality (\ref{nev_ukh_sum_abs_l_ij_greater2}) holds.
This yields
$|\det({\bf T})|<h_{100}$. 
Consequently, $|\det({\bf U})| < g_{101}$.
 For the corresponding simplex $S$, by \eqref{d_i_l_ij_formula},
$$
d_i(S) = 
\left\{
 \arraycolsep=3pt\def\arraystretch{1.2} 
\begin{array}{cl}
\frac{711}{719}=0.9888\ldots, & i=1,2,3, \\
\frac{237}{245}=0.9673\ldots, & i=4, \\
 1, & i=5,\ldots,100. \\
\end{array}
\right. 
$$
This means that $S$ is not a maximum volume simplex in  $Q_{100}$.


\section{ Problems for a Simplex and a Euclidean Ball}\label{nev_ukh_sec_9}

In this section, we will discuss the results of \cite{nevskii_mais_2018_25_6} for characteristics 
$\alpha(B_n;S)$ and $\xi(B_n;S)$ defined
for a simplex and the unit Euclidean ball in ${\mathbb R}^n$.
Replacing a cube with a ball makes many 
questions much more simpler. Nevertheless, geometric interpretation of general results
has a certain interest also in this particular case.
Moreother, we will note some additional applications of the basic Lagrange polynomials.

{\it The inradius of an $n$-dimensional simplex $S$} is the
maximum of the radii of balls contained within $S$.
The center of this unique maximum ball is called the {\it incenter of $S$.}
The boundary of the maximum ball is a sphere that has a single common point
with each $(n-1)$-dimensional face of $S$. By  {\it the circumradius of S} 
we mean the minimum of the radii of balls containing $S$.
The boundary of this unique minimal ball does not necessarily contain all the
vertices of $S$. Namely, this is only when the center of the minimal ball
lies inside the simplex.

The inradius $r$ and the circumradius $R$ of a simplex $S$ 
satisfy the so-called {\it Euler inequality}
\begin{equation}\label{euler_ineq}
R\geq nr.
\end{equation}
Equality in
(\ref{euler_ineq}) takes place if and only if 
$S$ is a regular simplex.
For history, proofs, and generalizations of Euler inequality, 
see, e.\,g.,
\cite{klamkin_1979},
\cite{yang_wang_1985},
\cite{vince_2008}.

  Let $x^{(1)},$ $\ldots,$
$x^{(n+1)}$ be the vertices and let $\lambda_1,$ $\ldots,$
$\lambda_{n+1}$ be the basic Lagrange polynomials of an nondegenerate simplex 
$S\subset {\mathbb R}^n$. Suppose
$\Gamma_j$ is an $(n-1)$-dimensional hyperplane given
by the equation $\lambda_j(x)=0$. By 
$\Sigma_j$ we mean an $(n-1)$-dimensional face of $S$ contained
in $\Gamma_j$. 
Symbol $h_j$ denotes the height of $S$ conducted from the vertex $x^{(j)}$ 
onto~$\Gamma_j$;
$r$ denotes the inradius of $S$. Define
$\sigma_j$ as the $(n-1)$-measure of $\Sigma_j$. Put
$\sigma:=\sum\limits_{j=1}^{n+1} \sigma_j$.
Consider the vector $a_j:=\{l_{1j},\ldots,l_{nj}\}$ orthogonal
to $\Gamma_j$ and directed into a subspace containing
$x^{(j)}$. Obviously,
$$\lambda_j(x)=
l_{1j}x_1+\ldots+
l_{nj}x_n+l_{n+1,j}=(a_j,x)+l_{n+1,j}=(a_j,x)+\lambda_j(0).$$

\smallskip
{\bf Theorem 10.1.}
{\it The following equalities are true:
\begin{equation}\label{alpha_bs_sum_l_ij_equality}
\alpha(B_n;S)=
\sum_{j=1}^{n+1}\left(\sum_{i=1}^n l_{ij}^2\right)^{1/2},
\end{equation}
\begin{equation}\label{alpha_bs_h_j_equality}
\alpha(B_n;S)=\sum_{j=1}^{n+1}\frac{1}{h_j},
\end{equation}
\begin{equation}\label{alpha_bs_1_r_equality}
\alpha(B_n;S)= \frac{1}{r},
\end{equation}
\begin{equation}\label{alpha_bs_sigma_nV}
\alpha(B_n;S)=\frac{\sigma}{n\vo(S)}.
\end{equation}
}

\smallskip
Let us note some corollaries of Theorem 10.1. We have
$$\frac{1}{r}=\sum_{j=1}^{n+1}\frac{1}{h_j}.$$
For proving, it is sufficient to apply
(\ref{alpha_bs_h_j_equality}) and
(\ref{alpha_bs_1_r_equality}). 
It seems to be interesting
that this relation (which evidently can be obtained also
in a direct way) occurs to be equivalent to the general formula \eqref{alpha_C_S_general_formula}  for
$\alpha(\Omega;S)$ in the 
case when a conveх body $\Omega$ coincide with the unit ball.

 The inradius $r$ and the incenter $z$ of  $S$ can be calculated
by formulae
\begin{equation}\label{r_formula}
r=\frac{1}{ \sum\limits_{j=1}^{n+1}\left(\sum\limits_{i=1}^n l_{ij}^2\right)^{1/2}},
\end{equation}
\begin{equation}\label{z_formula}
z=\frac{1}{ \sum\limits_{j=1}^{n+1}\left(\sum\limits_{i=1}^n l_{ij}^2\right)^{1/2}}
\sum\limits_{j=1}^{n+1}\left(\sum\limits_{i=1}^n l_{ij}^2\right)^{1/2} x^{(j)}.
\end{equation}
The tangent point of the ball $B(z;r)$ and the face
$\Sigma_k$ has the form
\begin{equation}\label{y_k_formula}
y^{(k)}=\frac{1}{ \sum\limits_{j=1}^{n+1}\left(\sum\limits_{i=1}^n l_{ij}^2\right)^{1/2}}
\left[\sum\limits_{j=1}^{n+1}\left(\sum\limits_{i=1}^n l_{ij}^2\right)^{1/2} x^{(j)}
-\frac{1}{\left(\sum\limits_{i=1}^n l_{ik}^2\right)^{1/2}} \left(l_{1k},\ldots,l_{nk}\right)
\right].
\end{equation}

Equality
(\ref{r_formula}) follows immediately from
(\ref{alpha_bs_sum_l_ij_equality}) and
(\ref{alpha_bs_1_r_equality}). To obtain
(\ref{z_formula}), let us remark that
$$r=
{\rm dist}(z;\Gamma_j)=
\frac{|\lambda_j(z)|}{\|a_j\|}.$$
Since $z$ lies inside $S$, each barycentric coordinate of this point 
$\lambda_j(z)$ is positive, i.\,e.,
$\lambda_j(z)=r\|a_j\|.$
Consequently,
$$z=\sum_{j=1}^{n+1}\lambda_j(z)x^{(j)}=
r\sum_{j=1}^{n+1} \|a_j\| x^{(j)}.$$
This coincides with (\ref{z_formula}). 
Finally, since
vector $a_k=\{l_{1k},\ldots,l_{nk}\}$
is orthogonal to 
$\Sigma_k$ and is directed
from this face inside the simplex, a unique common point of
$B(z;r)$ and $\Sigma_k$ has the form
$$y^{(k)}=z-\frac{r}{\|a_k\|}a_k=r\left( \sum_{j=1}^{n+1} \|a_j\| x^{(j)}-\frac{1}{\|a_k\|}
a_k\right).$$
The latter is equivalent to (\ref{y_k_formula}).

It is interesting to compare 
(\ref{alpha_bs_sum_l_ij_equality}) with  formula 
(\ref{alpha_q_prime_s_formula}) for $\alpha(Q_n^\prime;S)$. Since $B_n$ is~a~subset of 
$Q_n^\prime=[-1,1]^n$, we have $\alpha(B_n;S)\leq \alpha(Q_n^\prime;S)$.
Analytically, this inequality  also follows from the estimate
$$\left(\sum_{i=1}^n l_{ij}^2\right)^{1/2}\leq 
\sum_{i=1}^n |l_{ij}|.$$

For arbitrary $x^{(0)}$ and $\varrho>0$, the value
$\alpha\left(B(x^{(0)};\varrho);S\right)$ can be calculated
with the use of Theorem 10.1
and the equality 
$\alpha(B(x^{(0)};\varrho);S)$ $=$ $\varrho\alpha(B_n;S)$.

If $S\subset Q_n$, then $d_i(S)$ $\leq$ $1$.
As it was noted,
(\ref{alpha_d_i_formula}) immediately gives $\alpha(Q_n;S)\geq n$. 
The equality $\alpha(Q_n;S)=n$ holds when and only when each 
axial diameter of $S$ is equal to 1. 
The following proposition
expresses the analogues of these properties for simplices contained
in a ball.  This theorem is equivalent to the Euler inequality~(\ref{euler_ineq}).

\smallskip
{\bf Theorem 10.2.}
{\it If $S\subset B_n$, then $\alpha(B_n;S)\geq n.$ The equality
$\alpha(B_n;S)=n$ holds true if and only if 
$S$ is a regular simplex inscribed into $B_n$.
}

Now let us give a formula from \cite{nevskii_mais_2018_25_6} for the absorption index of a Euclidean ball by a~simplex.

\smallskip
{\bf Theorem 10.3.}
{\it Suppose $S$ is a nondegenerate simplex in
${\mathbb R}^n$, $x^{(0)}\in {\mathbb R}^n$,
$\varrho>0$. 
If $B\left(x^{(0)};\varrho\right)
 \not\subset S$, we have
\begin{equation}\label{ksi_b_x0_ro_s_l_ij_equality}
\xi\left(B\left(x^{(0)};\varrho\right);S\right)=
(n+1)\max_{1\leq j\leq n+1}
\left[\varrho\left(\sum_{i=1}^n l_{ij}^2\right)^{1/2}-
\sum_{i=1}^n l_{ij}x_i^{(0)}-l_{n+1,j}\right]+1.
\end{equation}
In particular, if
$B_n\not\subset S$, then
\begin{equation}\label{ksi_bs_l_ij_equality}
\xi(B_n;S)=
(n+1)\max_{1\leq j\leq n+1}\left[\left(\sum_{i=1}^n l_{ij}^2\right)^{1/2}
-l_{n+1,j}\right]+1.
\end{equation}
}

\smallskip
{\bf Theorem 10.4.}
{\it 
If $S\subset B_n$, then $\xi(B_n;S)\geq n.$ The equality
$\xi(B_n;S)=n$ takes place if and only if 
$S$ is a regular simplex inscribed into $B_n$.
}

\smallskip
Theorem 10.4 follows immediately from Theorem 10.2 
and the inequality
$\xi(B_n;S)\geq \alpha(B_n;S)$. Therefore,  $\xi_n(B_n)=n$. For proofs and commentaries, see \cite{nevskii_mais_2018_25_6}.

\section{Linear Interpolation on a Euclidean Ball} 
\label{nev_ukh_sec_11}

In the authors' paper \cite{nev_ukh_mais_2019_26_2}  and in the paper \cite{nevskii_mais_2019_26_3}
by M.\,Nevskii, some questions related to linear interpolation on a Euclidean ball in ${\mathbb R}^n$
are considered. Let us describe
the main results of these papers.

An interpolation projector $P:C(B_n)\to \Pi_1({\mathbb R}^n)$ is called {\it minimal} if
$\|P\|_{B_n}=\theta_n(B_n)$. The existence of a minimal projector follows from
the continuity of  $\|P\|_{B_n}$ as the nodes function which can be considered on a closed bounded subset of 
${\mathbb R}^m,$ $m=n(n+1)$ given by the conditions  $x^{(j)}\in B_n,$
$\det {\bf A}\geq \varepsilon_n>0$. It is showed that there exists a minimal projector with nodes belonging to the boundary sphere $\|x\|=1$. 

Suppose $B=B(x^{(0)};R)$, $P:C(B)\to \Pi_1({\mathbb R}^n)$ is an interpolation projector with nodes
 $x^{(j)}\in B,$ 
$\lambda_j(x)=
l_{1j}x_1+\ldots+
l_{nj}x_n+l_{n+1,j}$ are the basic polynomials of the simplex  $S$ having the vertices  $x^{(j)}$. 
Then
\begin{equation}\label{norm_P_B_Lagrange}
\|P\|_{B}=
\max_{x\in B}\sum_{j=1}^{n+1}
|\lambda_j(x)|,
\end{equation}
see (\ref{norm_P_intro_cepochka}) for $\Omega=B$. In \cite{nev_ukh_mais_2019_26_2} formula
(\ref{norm_P_B_Lagrange}) is supplemented with another expression for the projector norm.

 \smallskip 
{\bf Theorem 11.1.} 
 {\it   The following equality holds:
$$\|P\|_B=
\max\limits_{f_j=\pm 1} \left[
R \left(\sum_{i=1}^n\left(\sum_{j=1}^{n+1} f_jl_{ij}\right)^2\right)^{1/2}
+\left|\sum_{j=1}^{n+1}f_j\left(\sum_{i=1}^n l_{ij}x_i^{(0)}+l_{n+1,j}\right)\right|
\right]=$$
\begin{equation}\label{norm_P_B_formula}
=
\max\limits_{f_j=\pm 1} \left[
R \left(\sum_{i=1}^n\left(\sum_{j=1}^{n+1} f_jl_{ij}\right)^2\right)^{1/2}
+\left|\sum_{j=1}^{n+1}f_j\lambda_j(x^{(0)})\right|
\right].
\end{equation}
}

\smallskip
In the case when the center of gravity $c(S)$ of the simplex is coincided with the center of the ball, 
formula (\ref{norm_P_B_formula}) became noticeably simpler. Indeed,
if $c(S)=x^{(0)}$, then
$$\|P\|_B=
\max\limits_{f_j=\pm 1} \left[
R \left(\sum_{i=1}^n\left(\sum_{j=1}^{n+1} f_jl_{ij}\right)^2\right)^{1/2}
+\frac{1}{n+1}\left|\sum_{j=1}^{n+1}f_j\right|
\right].
$$

Let $S^*$ be a regular simplex inscribed into the ball
$B=$ $B(x^{(0)};R)$ and let
$P^*:C(B)\to \Pi_1({\mathbb R}^n)$ be the corresponding interpolation projector.
 It is easy to see that  $\|P^*\|_B$ 
  depends nor on the center $x^{(0)}$, nor on  
radius $R$ of the ball,  nor  on the choice of a regular simplex 
inscribed into that ball. 
In other words, $\|P^*\|_B$ is a function of  only dimension $n$. 
For $0\leq t\leq n+1$, consider the function 
$$\psi(t):=\frac{2\sqrt{n}}{n+1}\Bigl(t(n+1-t)\Bigr)^{1/2}+
\left|1-\frac{2t}{n+1}\right|.
$$
By definition, put
$a:=\left\lfloor\frac{n+1}{2}-\frac{\sqrt{n+1}}{2}\right\rfloor.$

 
  \smallskip 
{\bf Theorem 11.2.} {\it The following relations are true:
$$\|P^*\|_B=\max\{\psi(a),\psi(a+1)\},$$
\begin{equation}\label{norm_P_reg_ineqs}
\sqrt{n}\leq \|P^*\|_B\leq \sqrt{n+1}.
\end{equation}
Moreover,  $\|P^*\|_B=\sqrt{n}$ only for $n=1$;
the equality $\|P^*\|_B=\sqrt{n+1}$ holds if and only if 
$\sqrt{n+1}$ is an integer.} 

\smallskip
 

Recall that$\theta_n(B)$ denotes the minimal norm $\|P\|_B$ of an interpolation projector with nodes in $B$. Clearly,
$\theta_n(B)$ doesn't depend on the center and the radius of a ball; this is a function only of $n$. 
From \eqref{norm_P_reg_ineqs}, it follows immediately that $\theta_n(B)\geq c\sqrt{n}$.
As we will see, this estimate is exact in $n$. Further we take $B=B_n.$

With the use of Theorem 11.2  and arguments based on ~\eqref{nev_ksi_P_ineq} for  $\Omega=B_n$, it
was proved in \cite{nev_ukh_mais_2019_26_2} that
$$\theta_1(B_1)=1, \quad
\theta_2(B_2)=\frac{5}{3}, \quad
\theta_3(B_3)=2, \quad
\theta_4(B_4)=\frac{11}{5}.$$
For any $n\leq 4$, the equality $\|P\|_{B_n}=\theta_n(B_n)$ takes place only for projectors
corresponding to regular inscribed simplices.

Another approach to get  
$\|P^*\|_{B_n}=\theta_n(B_n)$, was suggested by M.\,Nevskii in \cite{nevskii_mais_2021_28_2}. This
approach is based on some geometric conjecture for a ball and an ellipsoid that are circumscribed around
a simplex.  However, also on this way, the equality   $\|P^*\|_{B_n}=\theta_n(B_n)$ 
so far is managed to be proved only for   $n\leq 4.$

 \begin{figure}[!htbp]
\begin{center}
\includegraphics[width=\textwidth]{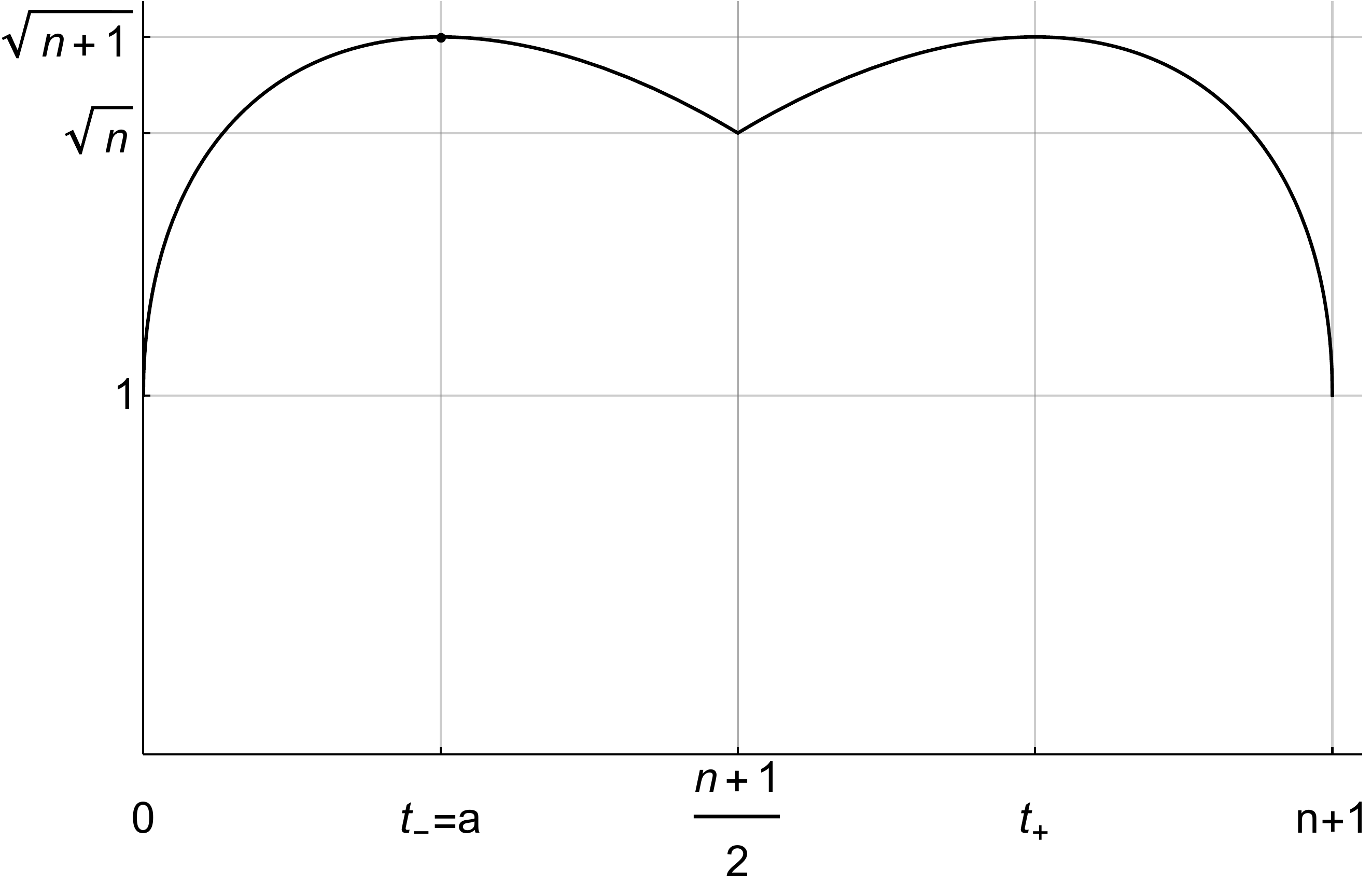}
\end{center}
\caption{The graph of $\psi(t)$ for $n=3$. Here $n+1=4$, $t_-=a=1$}
\label{fig:nev_ukl_n3}
\end{figure}

\begin{figure}[!htbp]
 \centering
\includegraphics[width=\textwidth]{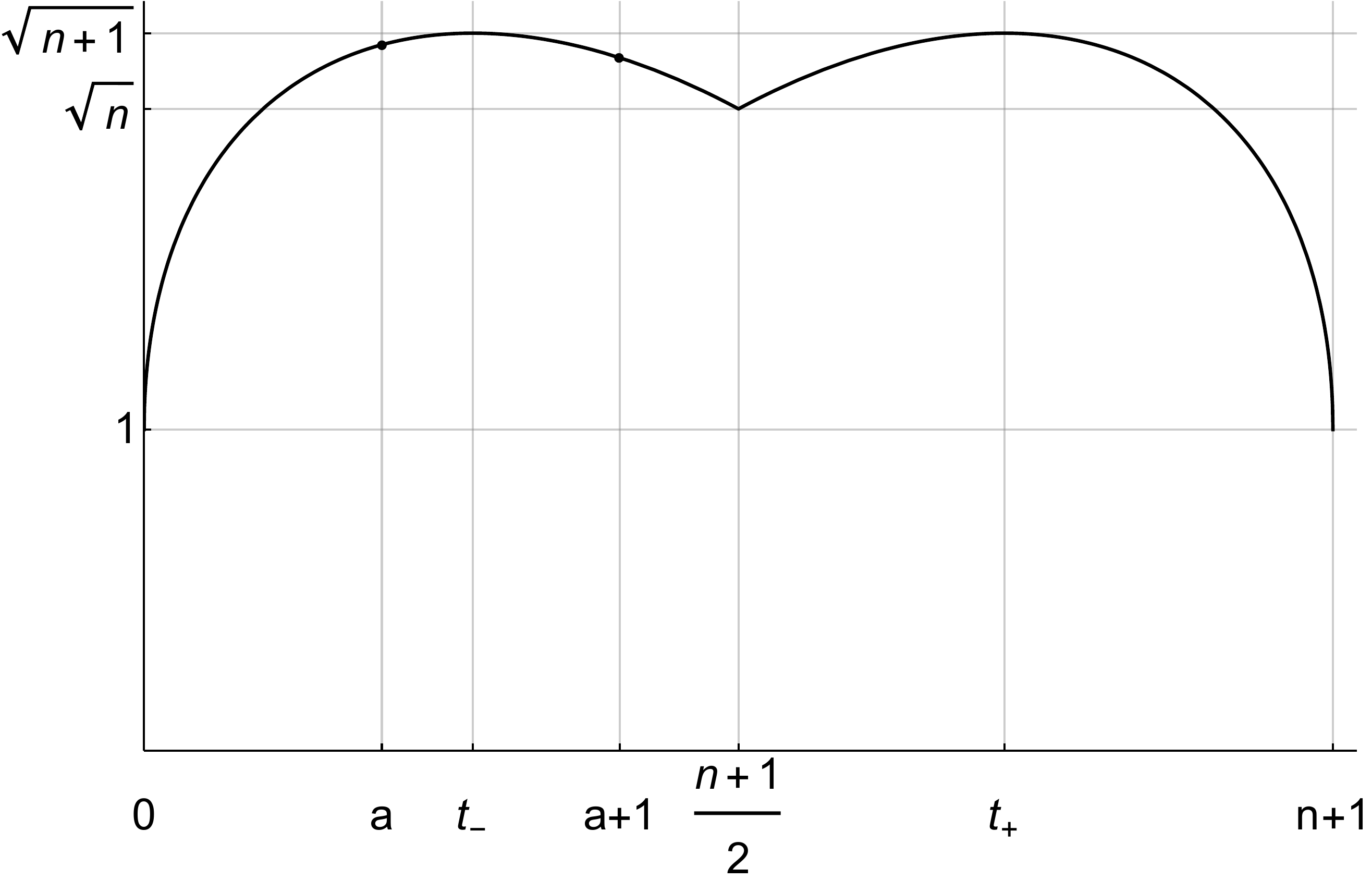}
\caption{The graph of  $\psi(t)$ for $n=4$. Here $n+1=5$, $a=1$, $t_-= \frac{ 5-\sqrt{5}}{2}$ 
}
\label{fig:nev_ukl_n4}
\end{figure}
\begin{figure}[!htbp]
 \centering
\includegraphics[width=\textwidth]{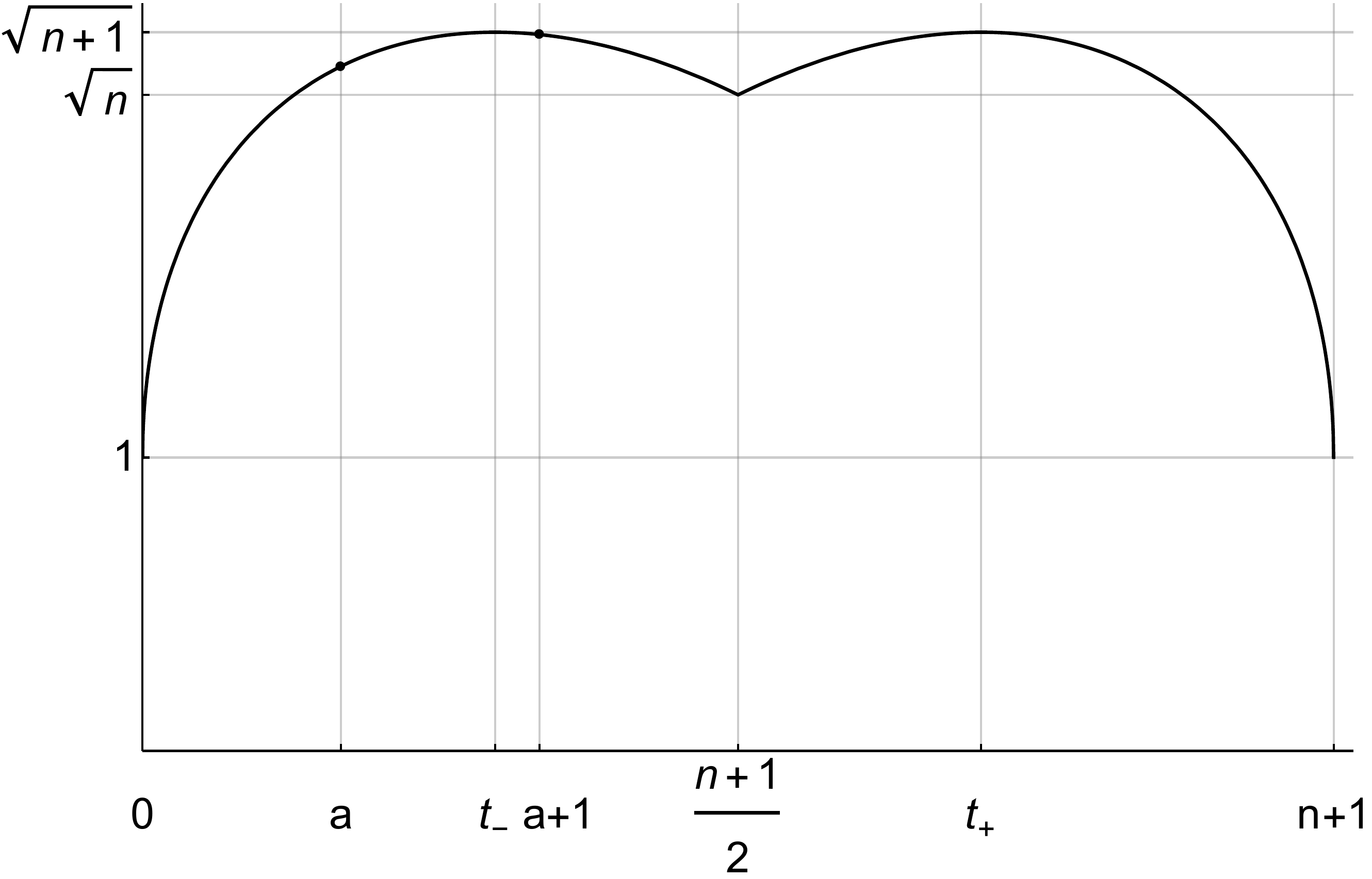}
\caption{The graph of $\psi(t)$ for$n=5$. Here $n+1=5$, $a=1$, $t_-= \frac{6-\sqrt{6}}{2}$} 
\label{fig:nev_ukl_n5}
\end{figure}
\begin{figure}[!htbp]
 \centering
\includegraphics[width=\textwidth]{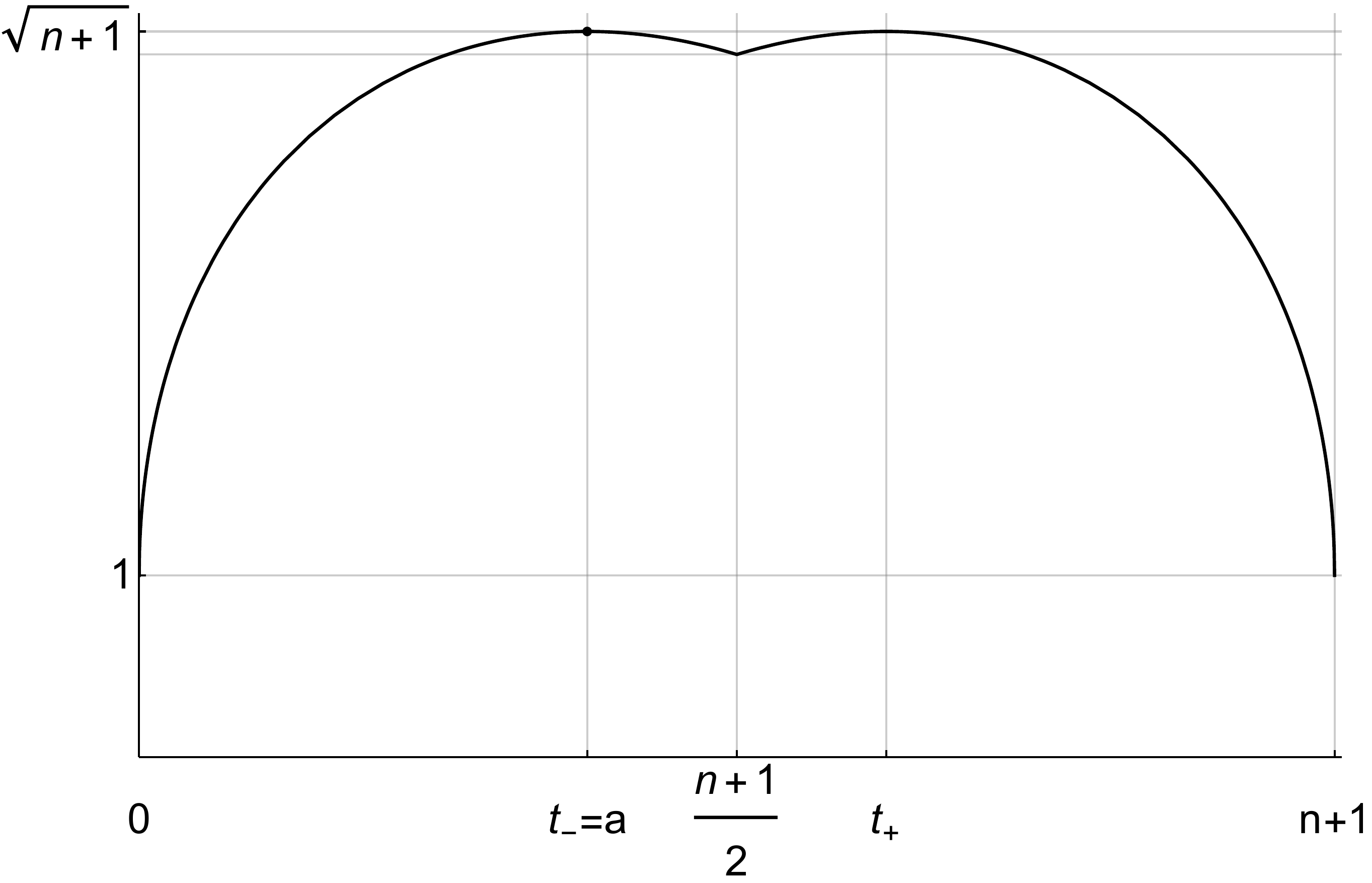}
\caption{The graph of $\psi(t)$ for $n=15$. Here $n+1=16$, $t_-=a=6$ }
\label{fig:nev_ukl_n15}
\end{figure}
%

Let us present some illustrations, the results of numerical analysis 
and some comments  given in  \cite{nev_ukh_mais_2019_26_2}. 
Graphs of the function 
$$ \psi(t)=\frac{2\sqrt{n}}{n+1}\Bigl(t(n+1-t)\Bigr)^{1/2}+
\left|1-\frac{2t}{n+1}\right|, \quad  0\leq t\leq n+1,$$
for $n=3, 4, 5, 15$ are shown in Fig.~\ref{fig:nev_ukl_n3}--\ref{fig:nev_ukl_n15}. We mark the maximum points
of this function $t_-:=\frac{n+1}{2}-\frac{\sqrt{n+1}}{2},$  
$t_+:=\frac{n+1}{2}+\frac{\sqrt{n+1}}{2}$ and also the points
$a=\left\lfloor\frac{n+1}{2}-
\frac{\sqrt{n+1}}{2}\right\rfloor$ and $a+1$. One  
of the points
 $a$ and $a+1$
maximizes 
$\psi(k)$ for integer $1\leq k\leq \frac{n+1}{2}$. 

\begin{figure}[!htbp]
 \centering
\includegraphics[width=\textwidth]{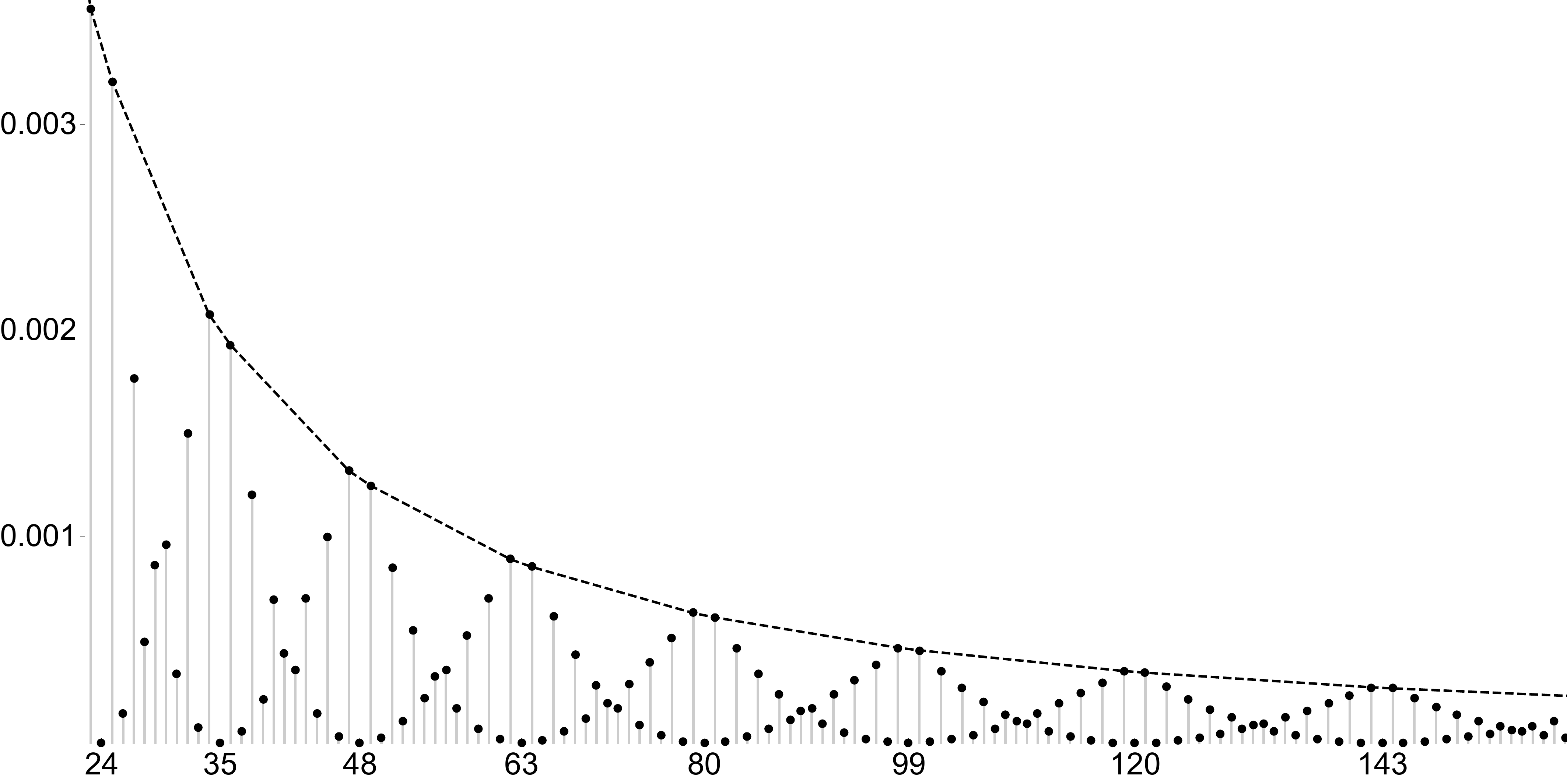}
\caption{The numbers $d_n=\sqrt{n+1}- \|P^*\|_{B_n}$ for $23\leq n\leq 160$}
\label{fig:nev_ukl_graph_diff}
\end{figure}

\begin{figure}[!htbp]
 \centering
\includegraphics[width=\textwidth]{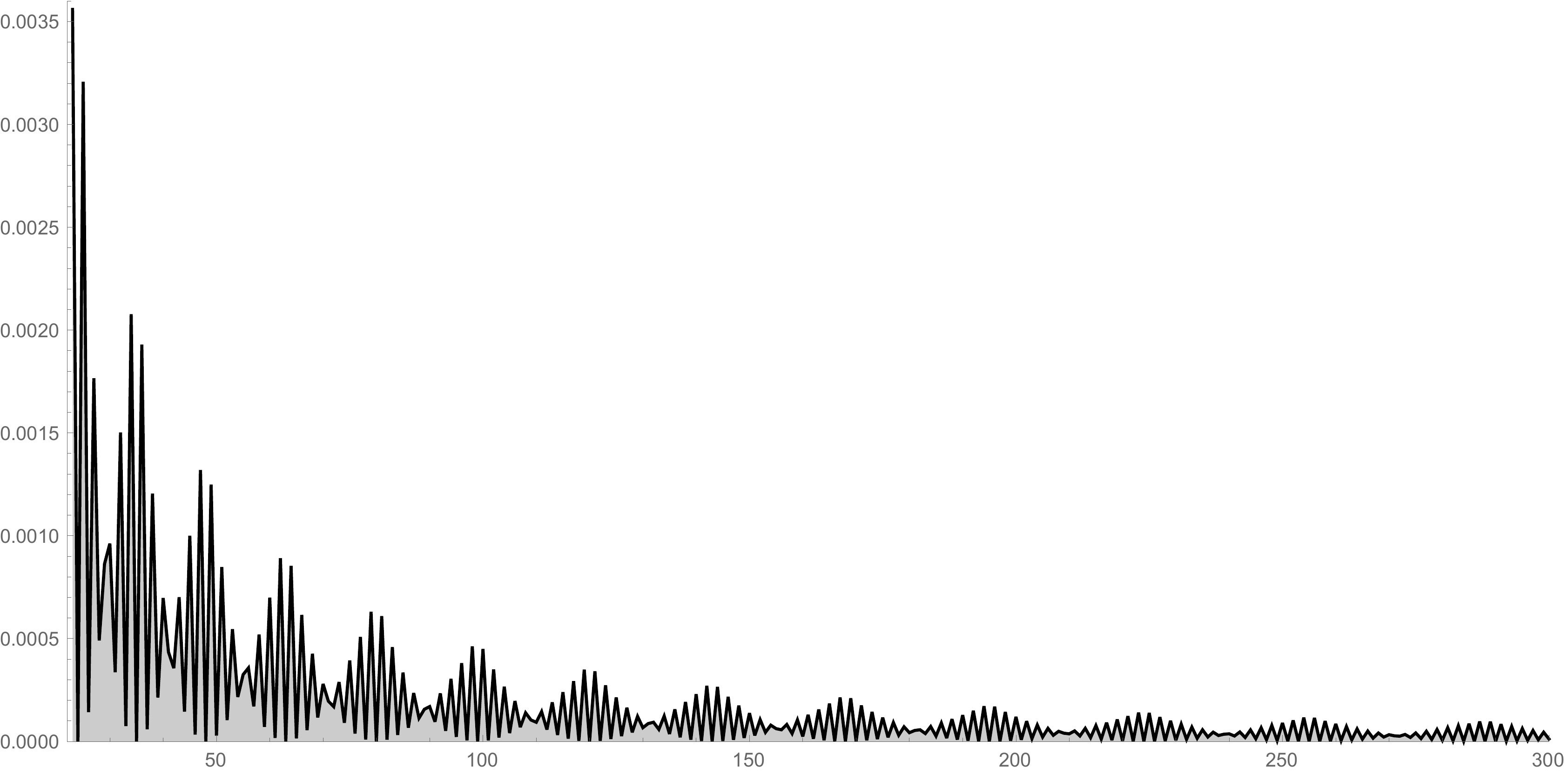}
\caption{The numbers $d_n=\sqrt{n+1}- \|P^*\|_{B_n}$ for $23\leq n\leq 300$}
\label{fig:nev_ukl_graph_diff_long}
\end{figure}

In Fig.~\ref{fig:nev_ukl_graph_diff}--\ref{fig:nev_ukl_graph_diff_long},
the values 
$d_n=\sqrt{n+1}-\|P^*\|_{B_n}$  for $n>23$ are presented.
As we have established, $\|P^*\|_{B_n}\leq\sqrt{n+1}$, with
an equality  only for $n=3, 8, 15, 24, 35, 48,$  $63, 80,$ $\ldots$, i.\,e.,
for $n$ having the form $ m ^ 2-1 $. 
It is at these points that $d_n=0$, as~can be seen in Fig.~\ref{fig:nev_ukl_graph_diff}.
The dashed line denotes the graph of the linear interpolation spline 
$l(n)$ constructed by the nodes $n=m^2-2,$ $n=m^2$
and by the values  $d_n$ in~these nodes.
Always $d_n\leq l(n)$; the equality comes only for $n=m^2-2$ and $n=m^2$.
We have $l(n)\to 0$ monotonically as $n\to\infty$.
Thus, the two-side estimate 
$$\sqrt{n+1}-l(n)\leq  \|P^*\|_{B_n} \leq \sqrt{n+1}$$
 takes place.
Both the right and left inequalities turn into equalities for an infinite set of~$n$.
In~\cite{nev_ukh_mais_2019_26_2}, we also give  the data related to the calculation of 
 $\|P^*\|_{B_n}$ 
\linebreak for $1\leq n\leq 15,$ $n=50, n=100,$ and~$n=1000$.

In \cite{nevskii_mais_2019_26_3}, M.\,Nevskii developed an approach wherein the norm of an interpolation
projector $P: C(B_n) \to \Pi_1({\mathbb R}^n)$ can be estimated from below
through the volume \linebreak of the corresponding simplex. 
The essential feature of this
approach is the application of the classical Legendre polynomials.
On this way, it was proved that $ \theta_n (B_n) \asymp \sqrt {n} $. 
Let us describe these results in more detail.

{\it The standardized Legendre polynomial of degree 
$n$} is the function
$$\chi_n(t):=\frac{1}{2^nn!}\left[ (t^2-1)^n \right] ^{(n)}.$$
The reasons why these polynomials appeared in the range of our questions are explained in Section 1.
There we also  gave the estimates for minimal norms of interpolation projectors
on a cube and on a ball with the participation of the function
$\chi_n^{-1}$ inverse to $\chi_n$ 
on the half-axis $ [1, + \ infty) $. Let us present estimates for interpolation on a ball proved 
в \cite{nevskii_mais_2019_26_3}.
The central result of this paper
 is the following theorem. Denote $\varkappa_n:=\vol(B_n)$.

\smallskip
{\bf Theorem 11.3.}
{\it Suppose
$P:C(B_n)\to\Pi_1({\mathbb R}^n)$ is an arbitrary
interpolation projector. Then for
the corresponding simplex 
$S\subset B_n$ and the vertices matrix ${\bf A}$ we have  
\begin{equation}\label{norm_P_vol_S_ineq}
\|P\|_{B_n}
\geq
\chi_n^{-1} 
\left(\frac{\varkappa_n}{\vo(S)}\right)=
\chi_n^{-1} 
\left(\frac{n!\varkappa_n}{|\det({\bf A})|}\right).
\end{equation}
}

\smallskip
Denote by
$\sigma_n$ the volume of a regular simplex inscribed into $B_n$. 
From~\eqref{norm_P_vol_S_ineq}, it follows  immediately   for each n
\begin{equation}\label{theta_n_chi_n_kappa_sigma_ineq}
\theta_n(B_n)
\geq
\chi_n^{-1} 
\left(\frac{\varkappa_n}{\sigma_n}\right).
\end{equation}
It is known that
$$\varkappa_n=\frac{\pi^
{\frac{n}{2}}}
{\Gamma\left(\frac{n}{2}+1\right)},\qquad
\sigma_n=\frac{1}{n!}\sqrt{n+1}\left(\frac{n+1}{n}\right)^{\frac{n}{2}},$$
$$\varkappa_{2k}=\frac{\pi^{k}}{k!},\qquad
\varkappa_{2k+1}=\frac{2^{k+1}\pi^{k}}{(2k+1)!!}=
\frac{2(k!)(4\pi)^k}{(2k+1)!}
$$
(see, e.\,g., \cite{fiht_2001},
\cite{nevskii_monograph}).
Therefore, the estimate (\ref{theta_n_chi_n_kappa_sigma_ineq})
can be made more concrete:
\begin{equation}\label{theta_n_chi_n_through_n_ineq}
\theta_n(B_n)
\geq
\chi_n^{-1} 
\left(\frac{\pi^{\frac{n}{2}}n!}{\Gamma\left(\frac{n}{2}+1\right)\sqrt{n+1}
\left(\frac{n+1}{n}\right)^{\frac{n}{2}}}
\right).
\end{equation}
If $n=2k$, then  (\ref{theta_n_chi_n_through_n_ineq})  is equivalent to the inequality
\begin{equation}\label{theta_n_chi_n_through_n_is_2k_ineq}
\theta_{2k}(B_{2k})
\geq
\chi_{2k}^{-1} 
\left(\frac{\pi^{k}(2k)!}{k!\sqrt{2k+1}
\left(\frac{2k+1}{2k}\right)^k}
\right).
\end{equation}
For $n=2k+1$, we have
\begin{equation}\label{theta_n_chi_n_through_n_is_2k_plus_1_ineq}
\theta_{2k+1}(B_{2k+1})
\geq
\chi_{2k+1}^{-1} 
\left(\frac{2(k!)(4\pi)^{k}}{\sqrt{2k+2}
\left(\frac{2k+2}{2k+1}\right)^{\frac{2k+1}{2}}}
\right).
\end{equation}

The  Stirling formula
$n!=\sqrt{2\pi n}\left(\frac{n}{e}\right)^n
e^{\frac{\zeta_n}{12n}}$,
$0<\zeta_n<1$, 
yields
\begin{equation}\label{n_fact_ineqs}
\sqrt{2\pi n}\left(\frac{n}{e}\right)^n<n!<\sqrt{2\pi n}\left(\frac{n}{e}\right)^n
e^{\frac{1}{12n}}.
\end{equation}
 It was proved in \cite[п.\,3.4.2]{nevskii_monograph} that 
\begin{equation}\label{chi_2k_2k_plus_1_ineqs}
\chi_{2k}^{-1}(s)>\left( \frac{(k!)^2s}
{ (2k)!}\right)^{\frac{1}{2k}}, \qquad
\chi_{2k+1}^{-1}(s)>\left( \frac{(k+1)!k!s}
{ (2k+1)!}\right)^{\frac{1}{2k+1}}. 
\end{equation}
With the use of the estimates \eqref{theta_n_chi_n_through_n_is_2k_ineq}--\eqref{chi_2k_2k_plus_1_ineqs}, the 
following result was obtained in   \cite{nevskii_mais_2019_26_3}.

\smallskip
{\bf Theorem 11.4.}
{\it There exists a constant
$c>0$ not depending on $n$ such that
\begin{equation}\label{theta_n_B_n_gt_c_sqrt_n}
\theta_n(B_n)>c\sqrt{n}.
\end{equation}
The inequality 
(\ref{theta_n_B_n_gt_c_sqrt_n}) takes place, e.\,g., with the constant
$$c=
\frac{  \sqrt[3]{\pi} }{\sqrt{12e}\cdot\sqrt[6]{3} } =  0.2135...$$
}

\smallskip
From Theorems 11.2 and 11.4, we obtain 
$\theta_n(B_n)\asymp \sqrt{n}$. 
Our results also  mean that any interpolation projector $P^*$
corresponding to an inscribed regular simplex 
has the norm equivalent to the minimal possible. In other words,
with constants not depending on $n$, we have
$\|P^*\|_{B_n}\asymp \theta_n(B_n)$.
The equality $\|P^*\|_{B_n}=\theta_n(B_n)$ remains proved
only for  $1\leq n\leq 4$.


\end{document}